\documentstyle[12pt,fleqn]{article}

\begin{document}

\title{Polar Complex Numbers in $n$ Dimensions}

\author{Silviu Olariu
\thanks{e-mail: olariu@ifin.nipne.ro}\\
Institute of Physics and Nuclear Engineering,\\
Department of Fundamental Experimental Physics\\
76900 Magurele, P.O. Box MG-6, Bucharest, Romania}

\date{4 August 2000}

\maketitle

\abstract

Polar commutative n-complex numbers of the form
$u=x_0+h_1x_1+h_2x_2+\cdots+h_{n-1}x_{n-1}$ are introduced in n dimensions, the
variables $x_0,...,x_{n-1}$ being real numbers.  The polar n-complex number can
be represented, in an even number of dimensions, by the modulus $d$, by the
amplitude $\rho$, by 2 polar angles $\theta_+,\theta_-$, by $n/2-2$ planar
angles $\psi_{k-1}$, and by $n/2-1$ azimuthal angles $\phi_k$. In an odd number
of dimensions, the polar n-complex number can be represented by $d, \rho$, by 1
polar angle $\theta_+$, by $(n-3)/2$ planar angles $\psi_{k-1}$, and by
$(n-1)/2$ azimuthal angles $\phi_k$.  The exponential function of a polar
n-complex number can be expanded in terms of the polar n-dimensional
cosexponential functions $g_{nk}(y), k=0,1,...,n-1$. Expressions are given for
these cosexponential functions.  The polar n-complex numbers can be written in
exponential and trigonometric forms with the aid of the modulus, amplitude and
the angular variables.  The polar n-complex functions defined by series of
powers are analytic, and the partial derivatives of the components of the
polar n-complex functions are closely related. The integrals of polar n-complex
functions are independent of path in regions where the functions are regular.
The fact that the exponential form of a polar n-complex numbers depends on the
cyclic variables $\phi_k$ leads to the concept of pole and residue for
integrals on closed paths.  The polynomials of polar n-complex variables can be
written as products of linear or quadratic factors, although the factorization
may not be unique.

\endabstract

\section{Introduction}

A regular, two-dimensional complex number $x+iy$ 
can be represented geometrically by the modulus $\rho=(x^2+y^2)^{1/2}$ and 
by the polar angle $\theta=\arctan(y/x)$. The modulus $\rho$ is multiplicative
and the polar angle $\theta$ is additive upon the multiplication of ordinary 
complex numbers.

The quaternions of Hamilton are a system of hypercomplex numbers
defined in four dimensions, the
multiplication being a noncommutative operation, \cite{1} 
and many other hypercomplex systems are
possible, \cite{2a}-\cite{2b} but these hypercomplex systems 
do not have all the required properties of regular, 
two-dimensional complex numbers which rendered possible the development of the 
theory of functions of a complex variable.

A system of complex numbers in $n$ dimensions is described in this work,
for which the multiplication is both associative and commutative, and which is 
rich enough in properties so that an exponential form exists and the concepts
of analytic n-complex 
function,  contour integration and residue can be defined.
The n-complex numbers introduced in this work have 
the form $u=x_0+h_1x_1+h_2x_2+\cdots+h_{n-1}x_{n-1}$, the variables 
$x_0,...,x_{n-1}$ being real numbers. The multiplication rules for the complex
units $h_1,...,h_{n-1}$ 
are $h_j h_k =h_{j+k}$ if $0\leq j+k\leq n-1$, and $h_jh_k=h_{j+k-n}$ if
$n\leq j+k\leq 2n-2$.
The product of two n-complex numbers is equal to zero if both numbers are
equal to zero, or if the numbers belong to certain n-dimensional hyperplanes
described further in this work. 

If the n-complex number $u=x_0+h_1x_1+h_2x_2+\cdots+h_{n-1}x_{n-1}$ is
represented by the point $A$ of coordinates $x_0,x_1,...,x_{n-1}$, 
the position of the point $A$ can be described, in an even number 
of dimensions, by the modulus $d=(x_0^2+x_1^2+\cdots+x_{n-1}^2)^{1/2}$, 
by $n/2-1$ azimuthal angles $\phi_k$, by $n/2-2$
planar angles $\psi_{k-1}$, and by 2 polar angles $\theta_+,\theta_-$. In an
odd number of dimensions, the position of the point $A$ is described by $d$, by
$(n-1)/2$ azimuthal angles $\phi_k$, by $(n-3)/2$ planar angles $\psi_{k-1}$,
and by 1 polar angle $\theta_+$. 
An amplitude $\rho$ can be defined for
even $n$ as $\rho^n=v_+v_-\rho_1^2\cdots\rho_{n/2-1}^2$,
and for odd $n$ as $\rho^n=v_+\rho_1^2\cdots\rho_{(n-1)/2}^2$,
where $v_+=x_0+\cdots+x_{n-1}, v_-=x_0-x_1+\cdots+x_{n-2}-x_{n-1}$,
and $\rho_k$ are radii in orthogonal two-dimensional planes defined further in
this work. The amplitude $\rho$, the variables $v_+, v_-$, the radii
$\rho_k$, the variables $(1/\sqrt{2})\tan\theta_+,(1/\sqrt{2})\tan\theta_-,
\tan\psi_{k-1}$ are multiplicative, and the azimuthal angles $\phi_k$ are
additive upon the multiplication of n-complex numbers.
Because of the role of the axis $v_+$ and, in
an even number of dimensions, of the axis $v_-$, in the description of the
position of the point $A$ with the aid of the polar angle $\theta_+$ and, in an
even number of dimensions, of the polar angle $\theta_-$, the hypercomplex
numbers studied in this work will be called polar n-complex number, to
distinguish them from the planar n-complex numbers, which exist in an even
number of dimensions. \cite{2u}

The exponential function of an n-complex number can be expanded in terms of
the polar n-dimensional cosexponential functions
$g_{nk}(y)=\sum_{p=0}^\infty y^{k+pn}/(k+pn)!$, $k=0,1,...,n-1$.
It is shown that
$g_{nk}(y) = \frac{1}{n}\sum_{l=0}^{n-1}$$\exp\left\{y\cos\left(2\pi l/n\right)
\right\} $$\cos\left\{y\sin\left(2\pi l/n\right)-2\pi kl/n\right\}$, 
$k=0,1,...,n-1$. Addition theorems and other relations
are obtained for the polar n-dimensional cosexponential functions.

The exponential form of an n-complex number, which in an even number of
dimensions $n$ can be defined for
$x_0+\cdots+x_{n-1}>0, x_0-x_1+\cdots+x_{n-2}-x_{n-1}>0$, is
$u=\rho \exp\left\{\sum_{p=1}^{n-1}h_p\left[
(1/n)\ln\sqrt{2}/\tan\theta_+
+((-1)^p/n)\ln\sqrt{2}/\tan\theta_-
-(2/n)\sum_{k=2}^{n/2-1}
\cos\left(2\pi kp/n\right)\ln\tan\psi_{k-1}
\right]\right\}$
 $\exp\left(\sum_{k=1}^{n/2-1}\tilde e_k\phi_k\right)$,
where $\tilde e_k=(2/n)\sum_{p=1}^{n-1}h_p\sin(2\pi pk/n)$. 
In an odd number of dimensions $n$, the exponential form exists for
$x_0+\cdots+x_{n-1}>0$, and is\\
$u=\rho \exp\left\{\sum_{p=1}^{n-1}h_p\left[
(1/n)\ln\sqrt{2}/\tan\theta_+
-(2/n)\sum_{k=2}^{(n-1)/2}
\cos\left(2\pi kp/n\right)\ln\tan\psi_{k-1}
\right]\right\}$
 $\exp\left(\sum_{k=1}^{(n-1)/2}\tilde e_k\phi_k\right)$. 
A trigonometric form also exists for an n-complex number $u$,
when $u$ is written as the product of the modulus $d$, of a factor depending on
the polar and planar angles $\theta_+, \theta_-, \psi_{k-1}$ and of an
exponential factor depending on the azimuthal angles $\phi_k$.

Expressions are given for the elementary functions of n-complex variable.
The functions $f(u)$ of n-complex variable which are defined by power series
have derivatives independent of the direction of approach to the point under
consideration. If the n-complex function $f(u)$ 
of the n-complex variable $u$ is written in terms of 
the real functions $P_k(x_0,...,x_{n-1}), k=0,...,n-1$, then
relations of equality  
exist between partial derivatives of the functions $P_k$. 
The integral $\int_A^B f(u) du$ of an n-complex
function between two points $A,B$ is independent of the path connecting $A,B$,
in regions where $f$ is regular.
If $f(u)$ is an analytic n-complex function, then 
$\oint_\Gamma f(u)du/(u-u_0)$
$=2\pi f(u_0)\sum_{k=1}^{[(n-1)/2]}\tilde e_k$ 
$\;{\rm int}(u_{0\xi_k\eta_k},\Gamma_{\xi_k\eta_k})$,
where the functional ${\rm int}$ takes the values 0 or 1 depending on the
relation between $u_{0\xi_k\eta_k}$ and $\Gamma_{\xi_k\eta_k}$, which are
respectively the projections of the point $u_0$ and of 
the loop $\Gamma$ on the plane defined by the orthogonal axes $\xi_k$ and
$\eta_k$, as expained further in this work.

A polynomial $u^m+a_1 u^{m-1}+\cdots+a_{m-1} u +a_m$ can be
written as a 
product of linear or quadratic factors, although the factorization may not be
unique. 

This paper belongs to a series of studies on commutative complex numbers in $n$
dimensions. \cite{2c} 
A detailed analysis of the cases for $n=2,3,4,5,6$  of the polar
n-complex numbers can be found in the corresponding
studies mentioned in Ref. \cite{2c}.

\section{Operations with polar n-complex numbers}

A complex number in $n$ dimensions is determined by its $n$ components
$(x_0,x_1,...,x_{n-1})$. The polar n-complex numbers and
their operations discussed in this work can be represented 
by  writing the n-complex number $(x_0,x_1,...,x_{n-1})$ as  
$u=x_0+h_1x_1+h_2x_2+\cdots+h_{n-1}x_{n-1}$, where 
$h_1, h_2, \cdots, h_{n-1}$ are bases for which the multiplication rules are 
\begin{equation}
h_j h_k =h_l ,\:l=j+k-n[(j+k)/n],
\label{1}
\end{equation}
for $ j,k,l=0,1,..., n-1$.
In Eq. (\ref{1}), $[(j+k)/n]$ denotes the integer part of $(j+k)/n$, 
the integer part being defined as $[a]\leq a<[a]+1$, so that
$0\leq j+k-n[(j+k)/n]\leq n-1$. 
In this work, brackets larger than the regular brackets
$[\;]$ do not have the meaning of integer part.
The significance of the composition laws in Eq.
(\ref{1}) can be understood by representing the bases $h_j, h_k$ by points on a
circle at the angles $\alpha_j=2\pi j/n,\alpha_k=2\pi k/n$, as shown in Fig. 1,
and the product $h_j h_k$ by the point of the circle at the angle 
$2\pi (j+k)/n$. If $2\pi\leq 2\pi (j+k)/n<4\pi$, the point represents the basis
$h_l$ of angle $\alpha_l=2\pi(j+k-n)/n$.

Two n-complex numbers 
$u=x_0+h_1x_1+h_2x_2+\cdots+h_{n-1}x_{n-1}$,
$u^\prime=x^\prime_0+h_1x^\prime_1+h_2x^\prime_2+\cdots+h_{n-1}x^\prime_{n-1}$ 
are equal if and only if $x_i=x^\prime_i, i=0,1,...,n-1$.
The sum of the n-complex numbers $u$
and
$u^\prime$ 
is
\begin{equation}
u+u^\prime=x_0+x^\prime_0+h_1(x_1+x^\prime_1)+\cdots
+h_{n-1}(x_{n-1} +x^\prime_{n-1}) .
\label{2}
\end{equation}
The product of the numbers $u, u^\prime$ is 
\begin{equation}
\begin{array}{l}
uu^\prime=x_0 x_0^\prime +x_1x_{n-1}^\prime
+x_2 x_{n-2}^\prime+x_3x_{n-3}^\prime
+\cdots+x_{n-1}x_1^\prime\\
+h_1(x_0 x_1^\prime+x_1x_0^\prime+x_2x_{n-1}^\prime+x_3x_{n-2}^\prime
+\cdots+x_{n-1} x_2^\prime) \\
+h_2(x_0 x_2^\prime+x_1x_1^\prime+x_2x_0^\prime+x_3x_{n-1}^\prime
+\cdots+x_{n-1} x_3^\prime) \\
\vdots\\
+h_{n-1}(x_0 x_{n-1}^\prime+x_1x_{n-2}^\prime
+x_2x_{n-3}^\prime+x_3x_{n-4}^\prime
+\cdots+x_{n-1} x_0^\prime).
\end{array}
\label{3}
\end{equation}
The product $uu^\prime$ can be written as
\begin{equation}
uu^\prime=\sum_{k=0}^{n-1}h_k\sum_{l=0}^{n-1}x_l x^\prime_{k-l+n[(n-k-1+l)/n]}.
\label{3a}
\end{equation}
If $u,u^\prime,u^{\prime\prime}$ are n-complex numbers, the multiplication 
is associative
\begin{equation}
(uu^\prime)u^{\prime\prime}=u(u^\prime u^{\prime\prime})
\label{3b}
\end{equation}
and commutative
\begin{equation}
u u^\prime=u^\prime u ,
\label{3c}
\end{equation}
because the product of the bases, defined in Eq. (\ref{1}), is associative and
commutative. The fact that the multiplication is commutative can be seen also
directly from Eq. (\ref{3}).
The n-complex zero is $0+h_1\cdot 0+\cdots+h_{n-1}\cdot 0,$ 
denoted simply 0, 
and the n-complex unity is $1+h_1\cdot 0+\cdots+h_{n-1}\cdot 0,$ 
denoted simply 1.

The inverse of the n-complex number $u=x_0+h_1x_1+h_2x_2+\cdots+h_{n-1}x_{n-1}$
is the n-complex number
$u^\prime
=x^\prime_0+h_1x^\prime_1+h_2x^\prime_2+\cdots+h_{n-1}x^\prime_{n-1}$ 
having the property that
\begin{equation}
uu^\prime=1 .
\label{4}
\end{equation}
Written on components, the condition, Eq. (\ref{4}), is
\begin{equation}
\begin{array}{l}
x_0 x_0^\prime +x_1x_{n-1}^\prime+x_2 x_{n-2}^\prime+x_3x_{n-3}^\prime
+\cdots+x_{n-1}x_1^\prime=1,\\
x_0 x_1^\prime+x_1x_0^\prime+x_2x_{n-1}^\prime+x_3x_{n-2}^\prime
+\cdots+x_{n-1} x_2^\prime=0,\\
x_0 x_2^\prime+x_1x_1^\prime+x_2x_0^\prime+x_3x_{n-1}^\prime
+\cdots+x_{n-1} x_3^\prime=0, \\
\vdots\\
x_0 x_{n-1}^\prime+x_1x_{n-2}^\prime+x_2x_{n-3}^\prime+x_3x_{n-4}^\prime
+\cdots+x_{n-1} x_0^\prime=0.
\end{array}
\label{5}
\end{equation}
The system (\ref{5}) has a solution provided that the determinant of the
system, 
\begin{equation}
\nu={\rm det}(A), 
\label{5b}
\end{equation}
is not equal to zero, $\nu\not=0$, where
\begin{equation}
A=\left(
\begin{array}{ccccc}
x_0     &   x_{n-1} &   x_{n-2}   & \cdots  &x_1\\
x_1     &   x_0     &   x_{n-1}   & \cdots  &x_2\\
x_2     &   x_1     &   x_0       & \cdots  &x_3\\
\vdots  &  \vdots   &  \vdots & \cdots  &\vdots \\
x_{n-1} &   x_{n-2} &   x_{n-3}   & \cdots  &x_0\\
\end{array}
\right).
\label{6}
\end{equation}
If $\nu>0$, the quantity 
\begin{equation}
\rho=\nu^{1/n}
\label{6b}
\end{equation}
will be called amplitude of the
n-complex number $u=x_0+h_1x_1+h_2x_2+\cdots+h_{n-1}x_{n-1}$.
The quantity $\nu$ can be written as a product of linear factors
\begin{equation}
\nu=\prod_{k=0}^{n-1}
\left(x_0+\epsilon_k x_1+\epsilon_k^2 x_2
+\cdots+\epsilon^{n-1}_k x_{n-1}\right),
\label{7}
\end{equation}
where $\epsilon_k=e^{2\pi ik/n}$, $i$ being the imaginary
unit. The factors appearing in Eq. (\ref{7}) are of the form 
\begin{equation}
x_0+\epsilon_k x_1+\epsilon_k^2 x_2+\cdots+\epsilon^{n-1}_k x_{n-1}
=v_k+i\tilde v_k,
\label{8}
\end{equation}
where
\begin{equation}
v_k=\sum_{p=0}^{n-1}x_p\cos\frac{2\pi kp}{n},
\label{9a}
\end{equation}
\begin{equation}
\tilde v_k=\sum_{p=0}^{n-1}x_p\sin\frac{2\pi kp}{n},
\label{9b}
\end{equation}
for $k=1,2,...,n-1$ and, if $n$ is even, $k\not=n/2$.
For $k=0$ the factor in Eq. (\ref{8}) is 
\begin{equation}
v_+=x_0+x_1+\cdots+x_{n-1},
\label{9bb}
\end{equation}
and if $n$ is even, for $k=n/2$ the factor in Eq. (\ref{8}) is
\begin{equation}
v_-=x_0-x_1+\cdots+x_{n-2}-x_{n-1}.
\label{9bbb}
\end{equation}
It can be seen that $v_k=v_{n-k}, \tilde v_k=-\tilde v_{n-k}$,
$k=1,...,[(n-1)/2]$. 
The variables $v_+, v_-, v_k, \tilde v_k, k=1,...,[(n-1)/2] $ 
will be called canonical polar n-complex variables. 
Therefore,  the factors appear in Eq. (\ref{7}) in complex-conjugate pairs of
the form $v_k+i\tilde v_k$ and $v_{n-k}+i\tilde v_{n-k}=v_k-i\tilde v_k$, where
$k=1,...,[(n-1)/2]$, so that the product $\nu$ is a real quantity. 
If $n$ is an even number, the quantity $\nu$ is
\begin{equation}
\nu=v_+v_-\prod_{k=1}^{n/2-1}(v_k^2+\tilde v_k^2), 
\label{9c}
\end{equation}
and if $n$ is an odd number, $\nu$ is
\begin{equation}
\nu=v_+\prod_{k=0}^{(n-1)/2}(v_k^2+\tilde v_k^2).
\label{9d}
\end{equation}
Thus, in an even number of dimensions $n$, an n-complex number has an inverse
unless it lies on one of the nodal hypersurfaces $x_0+x_1+\cdots+x_{n-1}=0$, or
$x_0-x_1+\cdots+x_{n-2}-x_{n-1}=0$, or $v_1=0, \tilde v_1=0$, ...,
or $v_{n/2-1}=0, \tilde v_{n/2-1}=0$. 
In an odd number of dimensions $n$, an n-complex number has an inverse
unless it lies on one of the nodal hypersurfaces $x_0+x_1+\cdots+x_{n-1}=0$,
or $v_1=0, \tilde v_1=0$, ..., or $v_{(n-1)/2}=0, \tilde v_{(n-1)/2}=0$. 

\section{Geometric representation of polar n-complex numbers}

The n-complex number $x_0+h_1x_1+h_2x_2+\cdots+h_{n-1}x_{n-1}$
can be represented by 
the point $A$ of coordinates $(x_0,x_1,...,x_{n-1})$. 
If $O$ is the origin of the n-dimensional space,  the
distance from the origin $O$ to the point $A$ of coordinates
$(x_0,x_1,...,x_{n-1})$ has the expression
\begin{equation}
d^2=x_0^2+x_1^2+\cdots+x_{n-1}^2 .
\label{10}
\end{equation}
The quantity $d$ will be called modulus of the n-complex number 
$u=x_0+h_1x_1+h_2x_2+\cdots+h_{n-1}x_{n-1}$. The modulus of an n-complex number
$u$ will be designated by $d=|u|$.

The exponential and trigonometric forms of the n-complex number $u$ can be
obtained conveniently in a rotated system of axes defined by a transformation
which, for even $n$, has the form
\begin{equation}
\left(
\begin{array}{c}
\xi_+\\
\xi_-\\
\vdots\\
\xi_k\\
\eta_k\\
\vdots
\end{array}\right)
=\left(
\begin{array}{ccccc}
\frac{1}{\sqrt{n}}&\frac{1}{\sqrt{n}}&
\cdots&\frac{1}{\sqrt{n}}&\frac{1}{\sqrt{n}}\\
\frac{1}{\sqrt{n}}&-\frac{1}{\sqrt{n}}&
\cdots&\frac{1}{\sqrt{n}}&-\frac{1}{\sqrt{n}}\\
\vdots&\vdots& &\vdots&\vdots\\
\sqrt{\frac{2}{n}}&\sqrt{\frac{2}{n}}
\cos\frac{2\pi k}{n}&\cdots&\sqrt{\frac{2}{n}}
\cos\frac{2\pi (n-2)k}{n}&\sqrt{\frac{2}{n}}\cos\frac{2\pi (n-1)k}{n}\\
0&\sqrt{\frac{2}{n}}\sin\frac{2\pi k}{n}&\cdots&
\sqrt{\frac{2}{n}}\sin\frac{2\pi (n-2)k}{n}&
\sqrt{\frac{2}{n}}\sin\frac{2\pi (n-1)k}{n}\\
\vdots&\vdots&&\vdots&\vdots
\end{array}
\right)
\left(
\begin{array}{c}
x_0\\x_1\\
\vdots\\ 
\vdots\\
\vdots\\
x_{n-1}
\end{array}
\right),
\label{11}
\end{equation}
where $k=1, 2, ... , n/2-1$.
For odd $n$ the rotation of the axes is described by the relations
\begin{equation}
\left(
\begin{array}{c}
\xi_+\\
\xi_1\\
\eta_1\\
\vdots\\
\xi_k\\
\eta_k\\
\vdots
\end{array}\right)
=\left(
\begin{array}{cccc}
\frac{1}{\sqrt{n}}&\frac{1}{\sqrt{n}}&\cdots&\frac{1}{\sqrt{n}}\\
\sqrt{\frac{2}{n}}&\sqrt{\frac{2}{n}}
\cos\frac{2\pi }{n}&\cdots&\sqrt{\frac{2}{n}}\cos\frac{2\pi (n-1)}{n}\\
0&\sqrt{\frac{2}{n}}\sin\frac{2\pi }{n}&\cdots&
\sqrt{\frac{2}{n}}\sin\frac{2\pi (n-1)}{n}\\
\vdots&\vdots& &\vdots\\
\sqrt{\frac{2}{n}}&\sqrt{\frac{2}{n}}
\cos\frac{2\pi k}{n}&\cdots&\sqrt{\frac{2}{n}}\cos\frac{2\pi (n-1)k}{n}\\
0&\sqrt{\frac{2}{n}}\sin\frac{2\pi k}{n}&\cdots&
\sqrt{\frac{2}{n}}\sin\frac{2\pi (n-1)k}{n}\\
\vdots&\vdots&&\vdots
\end{array}
\right)
\left(
\begin{array}{c}
x_0\\
x_1\\
x_2\\
\vdots\\ 
\vdots\\
\vdots\\
x_{n-1}
\end{array}
\right),
\label{12}
\end{equation}
where $k=0,1,...,(n-1)/2$.
The lines of the matrices in Eqs. (\ref{11}) or (\ref{12}) give the components
of the $n$ basis vectors of the new system of axes. These vectors have unit
length and are orthogonal to each other.
By comparing Eqs. (\ref{9a})-(\ref{9bbb}) and (\ref{11})-(\ref{12}) it can be
seen that
\begin{equation}
v_+= \sqrt{n}\xi_+ ,  v_-=  \sqrt{n}\xi_-,  
v_k= \sqrt{\frac{n}{2}}\xi_k , \tilde v_k= \sqrt{\frac{n}{2}}\eta_k ,
\label{12b}
\end{equation}
i.e. the two sets of variables differ only by scale factors.

The sum of the squares of the variables $v_k,\tilde v_k$ is, for even $n$,
\begin{equation}
\sum_{k=1}^{n/2-1}(v_k^2+\tilde v_k^2)=\frac{n-2}{2}(x_0^2+\cdots+x_{n-1}^2)
-2(x_0x_2+\cdots+x_{n-4}x_{n-2}+x_1x_3+\cdots+x_{n-3}x_{n-1}) ,
\label{13}
\end{equation}
and for odd $n$ the sum is
\begin{equation}
\sum_{k=1}^{(n-1)/2}(v_k^2+\tilde v_k^2)=\frac{n-1}{2}(x_0^2+\cdots+x_{n-1}^2)
-(x_0x_1+\cdots+x_{n-2}x_{n-1}) .
\label{14}
\end{equation}
The relation (\ref{13}) has been obtained with the aid of the identity, valid
for even $n$,
\begin{equation}
\sum_{k=1}^{n/2-1}\cos\frac{2\pi pk}{n}=\left\{
\begin{array}{l}
-1, \:\:{\rm for \:\;even}\:\:p,\\
0, \:\:{\rm for \:\;odd}\:\:p.
\end{array}
\right.
\label{15}
\end{equation}
The relation (\ref{14}) has been obtained with the aid of the identity, valid
for odd values of $n$,
\begin{equation}
\sum_{k=1}^{(n-1)/2}\cos\frac{2\pi pk}{n}=-\frac{1}{2}.
\label{16}
\end{equation}
From Eq. (\ref{13}) it results that, for even $n$,
\begin{equation}
d^2=\frac{1}{n}v_+^2+\frac{1}{n}v_-^2
+\frac{2}{n}\sum_{k=1}^{n/2-1}\rho_k^2,
\label{17}
\end{equation}
and from Eq. (\ref{14}) it results that, for odd $n$,
\begin{equation}
d^2=\frac{1}{n}v_+^2
+\frac{2}{n}\sum_{k=1}^{(n-1)/2}\rho_k^2.
\label{18}
\end{equation}
The relations (\ref{17}) and (\ref{18}) show that the square of the distance
$d$, Eq. (\ref{10}), is the sum of the squares of the projections
$v_+/\sqrt{n}, \rho_k\sqrt{2/n}$ and, for even $n$, of the square
of $v_-/\sqrt{n}$. This is consistent with the fact that the transformation in
Eqs. (\ref{11}) or (\ref{12}) is unitary.

The position of the point $A$ of coordinates $(x_0,x_1,...,x_{n-1})$ can be
also described with the aid of the distance $d$, Eq. (\ref{10}), and of $n-1$
angles defined further. Thus,
in the plane of the axes $v_k,\tilde v_k$, the radius $\rho_k$ and the
azimuthal angle 
$\phi_k$ can be introduced by the relations 
\begin{equation}
\rho_k^2=v_k^2+\tilde v_k^2, \:\cos\phi_k=v_k/\rho_k,\:\sin\phi_k=\tilde
v_k/\rho_k, 0\leq 
\phi_k<2\pi ,
\label{19a}
\end{equation}
so that there are $[(n-1)/2]$ azimuthal angles.
If the projection of the point $A$ on the plane of the axes $v_k,\tilde v_k$ is
$A_k$, 
and the projection of the point $A$ on the 4-dimensional space defined by the
axes $v_1, \tilde v_1, v_k,\tilde v_k$ is $A_{1k}$, the angle $\psi_{k-1}$
between the line 
$OA_{1k}$ and the 2-dimensional plane defined by the axes $v_k,\tilde v_k$ is 
\begin{equation}
\tan\psi_{k-1}=\rho_1/\rho_k, 
\label{19b}
\end{equation}
where $0\leq\psi_k\leq\pi/2, k=2,...,[(n-1)/2]$,
so that there are $[(n-3)/2]$ planar angles.
Moreover, there is a polar angle $\theta_+$, which can be
defined as the angle between the line $OA_{1+}$ and the axis $v_+$,
where $A_{1+}$ is the projection of the point $A$ on the 3-dimensional space
generated by the axes $v_1, \tilde v_1, v_+$,
\begin{equation}
\tan\theta_+=\frac{\sqrt{2}\rho_1}{v_+}, 
\label{19c}
\end{equation}
where $0\leq\theta_+\leq\pi$ , 
and in an even number of dimensions $n$ there is also a polar angle $\theta_-$,
which can be defined as the angle between the line $OA_{1-}$ and the axis
$v_-$, 
where $A_{1-}$ is the projection of the point $A$ on the 3-dimensional space
generated by the axes $v_1, \tilde v_1, v_-$,
\begin{equation}
\tan\theta_-=\frac{\sqrt{2}\rho_1}{v_-}, 
\label{19d}
\end{equation}
where $0\leq\theta_-\leq\pi$ .
In Eqs. (\ref{19c}) and (\ref{19d}), the factor $\sqrt{2}$ appears from the
ratio of the normalization factors in Eq. (\ref{12b}).
Thus, the position of the point $A$ is described, in an even number 
of dimensions, by the distance $d$, by $n/2-1$ azimuthal angles, by $n/2-2$
planar angles, and by 2 polar angles. In an odd number of dimensions, the
position of the point $A$ is described by $(n-1)/2$ azimuthal angles, by
$(n-3)/2$ planar angles, and by 1 polar angle. These angles are shown in Fig.
2. 

The variables $\rho_k$ can be expressed in terms of $d$ and the planar angles
$\psi_k$ as
\begin{equation}
\rho_k=\frac{\rho_1}{\tan\psi_{k-1}}, 
\label{20a}
\end{equation}
for $k=2,...,[(n-1)/2]$, where, for even $n$,
\begin{eqnarray}
\rho_1^2=\frac{nd^2}{2}
\left(\frac{1}{\tan^2\theta_+}+\frac{1}{\tan^2\theta_-}+1
+\frac{1}{\tan^2\psi_1}+\frac{1}{\tan^2\psi_2}+\cdots
+\frac{1}{\tan^2\psi_{n/2-2}}\right)^{-1},
\label{20b}
\end{eqnarray}
and for odd $n$
\begin{eqnarray}
\rho_1^2=\frac{nd^2}{2}
\left(\frac{1}{\tan^2\theta_+}+1
+\frac{1}{\tan^2\psi_1}+\frac{1}{\tan^2\psi_2}+\cdots
+\frac{1}{\tan^2\psi_{(n-3)/2}}\right)^{-1}.
\label{20c}
\end{eqnarray}

If
$u^\prime=x_0^\prime+h_1x_1^\prime+h_2x_2^\prime+\cdots+h_{n-1}x_{n-1}^\prime, 
u^{\prime\prime}=x^{\prime\prime}_0+h_1x^{\prime\prime}_1
+h_2x^{\prime\prime}_2+\cdots+h_{n-1}x^{\prime\prime}_{n-1}$ 
are n-complex numbers of parameters $v_+^\prime,  v_-^\prime,
\rho_k^\prime, \theta_+^\prime,  \theta_-^\prime,
\psi_k^\prime,\phi_k^\prime$ and respectively  
$v_+^{\prime\prime},  v_-^{\prime\prime},
\rho_k^{\prime\prime}, \theta_+^{\prime\prime}, \theta_-^\prime,
\psi_k^{\prime\prime}, \phi_k^{\prime\prime}$, then the parameters
$v_+, v_-, \rho_k,\theta_+,  \theta_-, \psi_k,\phi_k$ of the product n-complex
number $u=u^\prime u^{\prime\prime}$ are given by 
\begin{equation}
v_+=v_+^\prime v_+^{\prime\prime},  
\label{21a}
\end{equation}
\begin{equation}
\rho_k=\rho_k^\prime\rho_k^{\prime\prime}, 
\label{21b}
\end{equation}
for $k=1,..., [(n-1)/2]$,
\begin{equation}
\tan\theta_+=\frac{1}{\sqrt{2}}\tan\theta_+^\prime \tan\theta_+^{\prime\prime},
\label{21c}
\end{equation}
\begin{equation}
\tan\psi_k=\tan\psi_k^\prime \tan\psi_k^{\prime\prime},
\label{21d}
\end{equation}
for $k=1,...,[(n-3)/2]$,  
\begin{equation}
\phi_k=\phi_k^\prime+\phi_k^{\prime\prime}, 
\label{21e}
\end{equation}
for $k=1,...,[(n-1)/2]$, 
and, if $n$ is even,
\begin{equation}
v_-=v_-^\prime v_-^{\prime\prime},  
\label{21f}
\end{equation}
\begin{equation}
\tan\theta_-=\frac{1}{\sqrt{2}}\tan\theta_-^\prime \tan\theta_-^{\prime\prime}.
\label{21g}
\end{equation}
The Eqs. (\ref{21a}) and (\ref{21f}) can be checked directly, and  
Eqs. (\ref{21b})-(\ref{21e}) and (\ref{21g}) are a consequence of the relations
\begin{equation}
v_k=v_k^\prime v_k^{\prime\prime}
-\tilde v_k^\prime \tilde v_k^{\prime\prime},\;
\tilde v_k=v_k^\prime \tilde v_k^{\prime\prime}
+\tilde v_k^\prime v_k^{\prime\prime},
\label{22}
\end{equation}
and of the corresponding relations of definition. Then the product $\nu$ in
Eqs. (\ref{9c}) and (\ref{9d}) has the property that
\begin{equation}
\nu=\nu^\prime\nu^{\prime\prime} 
\label{23}
\end{equation}
and, if $\nu^\prime>0, \nu^{\prime\prime}>0$, the amplitude $\rho$ defined in
Eq. (\ref{6b}) has the property that
\begin{equation}
\rho=\rho^\prime\rho^{\prime\prime} .
\label{24}
\end{equation}

The fact that the amplitude of the product is equal to the product of the 
amplitudes, as written in Eq. (\ref{24}), can 
be demonstrated also by using a representation of the 
n-complex numbers by matrices, in which the n-complex number 
$u=x_0+h_1x_1+h_2x_2+\cdots+h_{n-1}x_{n-1}$ is represented by the matrix
\begin{equation}
U=\left(
\begin{array}{ccccc}
x_0     &   x_1     &   x_2   & \cdots  &x_{n-1}\\
x_{n-1} &   x_0     &   x_1   & \cdots  &x_{n-2}\\
x_{n-2} &   x_{n-1} &   x_0   & \cdots  &x_{n-3}\\
\vdots  &  \vdots   &  \vdots & \cdots  &\vdots \\
x_1     &   x_2     &   x_3   & \cdots  &x_0\\
\end{array}
\right).
\label{24b}
\end{equation}
The product $u=u^\prime u^{\prime\prime}$ is
represented by the matrix multiplication $U=U^\prime U^{\prime\prime}$.
The relation (\ref{23}) is then a consequence of the fact the determinant 
of the product of matrices is equal to the product of the determinants 
of the factor matrices. The use of the representation of the n-complex
numbers with matrices provides an alternative demonstration of the fact
that the product of n-complex numbers is associative, as
stated in Eq. (\ref{3b}).

According to Eqs. (\ref{21a}, (\ref{21b}), (\ref{21f}), (\ref{17}) and
(\ref{18}), the modulus of the product $uu^\prime$ is, for even $n$, 
\begin{equation}
|uu^\prime|^2=
\frac{1}{n}(v_+v_+^\prime)^2+\frac{1}{n}(v_-v_-^\prime)^2
+\frac{2}{n}\sum_{k=1}^{n/2-1}(\rho_k\rho_k^\prime)^2 ,
\label{25a}
\end{equation}
and for odd $n$
\begin{equation}
|uu^\prime|^2=
\frac{1}{n}(v_+v_+^\prime)^2
+\frac{2}{n}\sum_{k=1}^{(n-1)/2}(\rho_k\rho_k^\prime)^2 .
\label{25b}
\end{equation}
Thus, if the product of two n-complex numbers is zero, $uu^\prime=0$, then
$v_+v_+^\prime=0, \rho_k\rho_k^\prime=0, k=1,...,[(n-1)/2]$ and, if $n$ is
even, $v_-v_-^\prime=0$. This means that either $u=0$, or $u^\prime=0$, or $u,
u^\prime$ belong to orthogonal hypersurfaces in such a way that the
afore-mentioned products of components should be equal to zero.

\section{The polar n-dimensional cosexponential functions}

The exponential function of the n-complex variable $u$ can be defined by the
series
\begin{equation}
\exp u = 1+u+u^2/2!+u^3/3!+\cdots . 
\label{26}
\end{equation}
It can be checked by direct multiplication of the series that
\begin{equation}
\exp(u+u^\prime)=\exp u \cdot \exp u^\prime . 
\label{27}
\end{equation}
If $u=x_0+h_1x_1+h_2x_2+\cdots+h_{n-1}x_{n-1}$,
then $\exp u$ can be calculated as 
$\exp u=\exp x_0 \cdot \exp (h_1x_1) \cdots \exp (h_{n-1}x_{n-1})$.

It can be seen with the aid of the representation in Fig. 1 that 
\begin{equation}
h_k^{n+p}=h_k^p, \:p\:\:{\rm integer},
\label{28}
\end{equation}
for $ k=1,...,n-1$.
Then $e^{h_k y}$ can be written as
\begin{equation}
e^{h_k y}=\sum_{p=0}^{n-1}h_{kp-n[kp/n]}g_{np}(y),
\label{28b}
\end{equation}
where the expression of the functions $g_{nk}$, which will be
called polar cosexponential functions in $n$ dimensions, is
\begin{equation}
g_{nk}(y)=\sum_{p=0}^\infty y^{k+pn}/(k+pn)!, 
\label{29}
\end{equation}
for $ k=0,1,...,n-1$.

If $n$ is even, the polar cosexponential functions of even index $k$ are
even functions, $g_{n,2p}(-y)=g_{n,2p}(y)$, $p=0,1,...,n/2-1$,
and the polar cosexponential functions of odd index 
are odd functions, $g_{n,2p+1}(-y)=-g_{n,2p+1}(y)$, $p=0,1,...,n/2-1$. For odd
values of $n$, the polar cosexponential functions do not have a definite
parity. It can be checked that
\begin{equation}
\sum_{k=0}^{n-1}g_{nk}(y)=e^y
\label{29a}
\end{equation}
and, for even $n$, 
\begin{equation}
\sum_{k=0}^{n-1}(-1)^k g_{nk}(y)=e^{-y}.
\label{29b}
\end{equation}

The expression of the polar n-dimensional cosexponential functions is
\begin{equation}
g_{nk}(y)=\frac{1}{n}\sum_{l=0}^{n-1}
\exp\left[y\cos\left(\frac{2\pi l}{n}\right)
\right]
\cos\left[y\sin\left(\frac{2\pi l}{n}\right)-\frac{2\pi kl}{n}\right], 
\label{30}
\end{equation}
for $k=0,1,...,n-1$.
In order to check that the function in Eq. (\ref{30}) has the series expansion
written in Eq. (\ref{29}), the right-hand side of Eq. (\ref{30}) will be
written as
\begin{equation}
g_{nk}(y)=\frac{1}{n}\sum_{l=0}^{n-1}{\rm Re}\left\{
\exp\left[\left(\cos\frac{2\pi l}{n}+i\sin\frac{2\pi l}{n}\right)y
-i\frac{2\pi kl}{n}\right]\right\}, 
\label{31}
\end{equation}
for $k=0,1,...,n-1$,
where ${\rm Re} (a+ib)=a$, with $a$ and $b$ real numbers. The part of the
exponential depending on $y$ can be expanded in a series,
\begin{equation}
g_{nk}(y)=\frac{1}{n}\sum_{p=0}^\infty\sum_{l=0}^{n-1}{\rm Re}
\left\{\frac{1}{p!}
\exp\left[i\frac{2\pi l}{n}(p-k)\right]y^p\right\}, 
\label{32}
\end{equation}
for $k=0,1,...,n-1$.
The expression of $g_{nk}(y)$ becomes
\begin{equation}
g_{nk}(y)=\frac{1}{n}\sum_{p=0}^\infty\sum_{l=0}^{n-1}
\left\{\frac{1}{p!}
\cos\left[\frac{2\pi l}{n}(p-k)\right]y^p\right\}, 
\label{33}
\end{equation}
for $k=0,1,...,n-1$ and, since
\begin{equation}
\frac{1}{n}\sum_{l=0}^{n-1}\cos\frac{2\pi l}{n}(p-k)
=\left\{
\begin{array}{l}
1, \:\:{\rm if}\:\: p-k\:\: {\rm is \:\;a\:\;multiple\:\;of\:\;}n,\\
0, \:\:{\rm otherwise},
\end{array}
\right.
\label{34}
\end{equation}
this yields indeed the expansion in Eq. (\ref{29}).

It can be shown from Eq. (\ref{30}) that
\begin{equation}
\sum_{k=0}^{n-1}g_{nk}^2(y)=\frac{1}{n}\sum_{l=0}^{n-1}\exp\left[2y\cos\left(
\frac{2\pi l}{n}\right)\right].
\label{34a}
\end{equation}
It can be seen that the right-hand side of Eq. (\ref{34a}) does not contain
oscillatory terms. If $n$ is a multiple of 4, it can be shown by replacing $y$
by $iy$ in Eq. (\ref{34a}) that
\begin{equation}
\sum_{k=0}^{n-1}(-1)^kg_{nk}^2(y)=\frac{2}{n}\left\{1+\cos 2y
+\sum_{l=1}^{n/4-1}\cos\left[2y\cos\left(\frac{2\pi l}{n}
\right)\right]\right\},
\label{34b}
\end{equation}
which does not contain exponential terms.

Addition theorems for the polar n-dimensional cosexponential functions can be
obtained from the relation $\exp h_1(y+z)=\exp h_1 y \cdot\exp h_1 z $, by
substituting the expression of the exponentials as given in Eq. (\ref{28b})
for $k=1$, $e^{h_1 y}=g_{n0}(y)+h_1g_{n1}(y)+\cdots+h_{n-1} g_{n,n-1}(y)$,
\begin{eqnarray}
\lefteqn{g_{nk}(y+z)=g_{n0}(y)g_{nk}(z)
+g_{n1}(y)g_{n,k-1}(z)+\cdots+g_{nk}(y)g_{n0}(z)\nonumber}\\
&&+g_{n,k+1}(y)g_{n,n-1}(z)+g_{n,k+2}(y)g_{n,n-2}(z)
+\cdots+g_{n,n-1}(y)g_{n,k+1}(z) ,
\label{35a}
\end{eqnarray}
where $k=0,1,...,n-1$.
For $y=z$ the relations (\ref{35a}) take the form
\begin{eqnarray}
\lefteqn{g_{nk}(2y)=g_{n0}(y)g_{nk}(y)
+g_{n1}(y)g_{n,k-1}(y)+\cdots+g_{nk}(y)g_{n0}(y)\nonumber}\\
&&+g_{n,k+1}(y)g_{n,n-1}(y)+g_{n,k+2}(y)g_{n,n-2}(y)
+\cdots+g_{n,n-1}(y)g_{n,k+1}(y) ,
\label{35b}
\end{eqnarray}
where $k=0,1,...,n-1$.
For $y=-z$ the relations (\ref{35a}) and (\ref{29}) yield
\begin{equation}
g_{n0}(y)g_{n0}(-y)
+g_{n1}(y)g_{n,n-1}(-y)+g_{n2}(y)g_{n,n-2}(-y)
+\cdots+g_{n,n-1}(y)g_{n1}(-y)=1 ,
\label{36a}
\end{equation}
\begin{eqnarray}
\lefteqn{g_{n0}(y)g_{nk}(-y)+g_{n1}(y)g_{n,k-1}(-y)
+\cdots+g_{nk}(y)g_{n0}(-y)\nonumber}\\
&&+g_{n,k+1}(y)g_{n,n-1}(-y)+g_{n,k+2}(y)g_{n,n-2}(-y)
+\cdots+g_{n,n-1}(y)g_{n,k+1}(-y)=0 ,
\nonumber\\
&&
\label{36b}
\end{eqnarray}
for $k=1,...,n-1$.

From Eq. (\ref{28b}) it can be shown, for natural numbers $l$, that
\begin{equation}
\left(\sum_{p=0}^{n-1}h_{kp-n[kp/n]}g_{np}(y)\right)^l
=\sum_{p=0}^{n-1}h_{kp-n[kp/n]}g_{np}(ly), 
\label{37}
\end{equation}
where $k=0,1,...,n-1$.
For $k=1$ the relation (\ref{37}) is
\begin{equation}
\left\{g_{n0}(y)+h_1g_{n1}(y)+\cdots+h_{n-1}g_{n,n-1}(y)\right\}^l
=g_{n0}(ly)+h_1g_{n1}(ly)+\cdots+h_{n,n-1}g_{n,n-1}(ly).
\label{37b}
\end{equation}

If
\begin{equation}
a_k=\sum_{p=0}^{n-1}g_{np}(y)\cos\left(\frac{2\pi kp}{n}\right), 
\label{38a}
\end{equation}
for $k=0,1,...,n-1$, and
\begin{equation}
b_k=\sum_{p=0}^{n-1}g_{np}(y)\sin\left(\frac{2\pi kp}{n}\right), 
\label{38b}
\end{equation}
for $k=1,...,n-1$, 
where $g_{nk}(y)$ are the polar cosexponential functions in Eq. (\ref{30}), it
can 
be shown that 
\begin{equation}
a_k=\exp\left[y\cos\left(\frac{2\pi k}{n}\right)\right]
\cos\left[y\sin\left(\frac{2\pi k}{n}\right)\right], 
\label{39a}
\end{equation}
where $k=0,1,...,n-1$,
\begin{equation}
b_k=\exp\left[y\cos\left(\frac{2\pi k}{n}\right)\right]
\sin\left[y\sin\left(\frac{2\pi k}{n}\right)\right], 
\label{39b}
\end{equation}
where $k=1,...,n-1$.
If
\begin{equation}
G_k^2=a_k^2+b_k^2, 
\label{40}
\end{equation}
for $k=1,...,n-1$,
then from Eqs. (\ref{39a}) and (\ref{39b}) it results that
\begin{equation}
G_k^2=\exp\left[2y\cos\left(\frac{2\pi k}{n}\right)\right], 
\label{41}
\end{equation}
where $k=1,...,n-1$. If
\begin{equation}
G_+=g_{n0}+g_{n1}+\cdots+g_{n,n-1} ,
\label{42a}
\end{equation}
from Eq. (\ref{38a}) it results that $G_+=a_0$, so that $G_+=e^y$, and, in an
even number of dimensions $n$, if 
\begin{equation}
G_-=g_{n0}-g_{n1}+\cdots+g_{n,n-2}-g_{n,n-1} ,
\label{42b}
\end{equation}
from Eq. (\ref{38a}) it results that $G_-=a_{n/2}$, so that $G_{n/2}=e^{-y}$.
Then with the aid of Eq. (\ref{15}) applied for $p=1$ it can be shown that the
polar n-dimensional cosexponential functions have the property that, 
for even $n$,
\begin{equation}
G_+G_-\prod_{k=1}^{n/2-1}G_k^2=1,
\label{43a}
\end{equation}
and in an odd number of dimensions, with the aid of Eq. (\ref{16}) it can be
shown that  
\begin{equation}
G_+\prod_{k=1}^{(n-1)/2}G_k^2=1.
\label{43b}
\end{equation}

The polar n-dimensional cosexponential functions are solutions of the
$n^{\rm th}$-order differential equation
\begin{equation}
\frac{d^n\zeta}{du^n}=\zeta ,
\label{44}
\end{equation}
whose solutions are of the form
$\zeta(u)=A_0g_{n0}(u)+A_1g_{n1}(u)+\cdots+A_{n-1}g_{n,n-1}(u).$ 
It can be checked that the derivatives of the polar cosexponential functions
are related by
\begin{equation}
\frac{dg_{n0}}{du}=g_{n,n-1}, \:
\frac{dg_{n1}}{du}=g_{n0}, \:...,
\frac{dg_{n,n-2}}{du}=g_{n,n-3} ,
\frac{dg_{n,n-1}}{du}=g_{n,n-2} .
\label{45}
\end{equation}

\section{Exponential and trigonometric forms of polar n-complex numbers}

In order to obtain the exponential and trigonometric forms of n-complex
numbers, a canonical base 
$e_+,e_-,e_1,\tilde e_1,...,e_{n/2-1},\tilde e_{n/2-1}$ 
for the polar n-complex numbers will be introduced for
even $n$ by the relations 
\begin{equation}
\left(
\begin{array}{c}
e_+\\
e_-\\
\vdots\\
e_k\\
\tilde e_k\\
\vdots
\end{array}\right)
=\left(
\begin{array}{ccccc}
\frac{1}{n}&\frac{1}{n}&\cdots&\frac{1}{n}&\frac{1}{n}\\
\frac{1}{n}&-\frac{1}{n}&\cdots&\frac{1}{n}&-\frac{1}{n}\\
\vdots&\vdots& &\vdots&\vdots\\
\frac{2}{n}&\frac{2}{n}\cos\frac{2\pi k}{n}&\cdots&
\frac{2}{n}\cos\frac{2\pi (n-2)k}{n}&\frac{2}{n}\cos\frac{2\pi (n-1)k}{n}\\
0&\frac{2}{n}\sin\frac{2\pi k}{n}&\cdots&
\frac{2}{n}\sin\frac{2\pi (n-2)k}{n}&\frac{2}{n}\sin\frac{2\pi (n-1)k}{n}\\
\vdots&\vdots&&\vdots&\vdots
\end{array}
\right)
\left(
\begin{array}{c}
1\\h_1\\
\vdots\\ 
\vdots\\
\vdots\\
h_{n-1}
\end{array}
\right),
\label{e11}
\end{equation}
where $k=1, 2, ... , n/2-1$.
For odd $n$, the canonical base 
$e_+, e_1,\tilde e_1,...e_{(n-1)/2},\tilde e_{(n-1)/2}$ 
for the polar n-complex numbers will be introduced by the relations
\begin{equation}
\left(
\begin{array}{c}
e_+\\
e_1\\
\tilde e_1\\
\vdots\\
e_k\\
\tilde e_k\\
\vdots
\end{array}\right)
=\left(
\begin{array}{ccccc}
\frac{1}{n}&\frac{1}{n}&\cdots&\frac{1}{n}\\
\frac{2}{n}&\frac{2}{n}\cos\frac{2\pi }{n}&\cdots&
\frac{2}{n}\cos\frac{2\pi (n-1)}{n}\\
0&\frac{2}{n}\sin\frac{2\pi }{n}&\cdots&\frac{2}{n}\sin\frac{2\pi (n-1)}{n}\\
\vdots&\vdots& &\vdots\\
\frac{2}{n}&\frac{2}{n}\cos\frac{2\pi k}{n}&
\cdots&\frac{2}{n}\cos\frac{2\pi (n-1)k}{n}\\
0&\frac{2}{n}\sin\frac{2\pi k}{n}&\cdots&\frac{2}{n}\sin\frac{2\pi (n-1)k}{n}\\
\vdots&\vdots&&\vdots
\end{array}
\right)
\left(
\begin{array}{c}
1\\
h_1\\
h_2\\
\vdots\\ 
\vdots\\
\vdots\\
h_{n-1}
\end{array}
\right),
\label{e12}
\end{equation}
where $k=0,1,...,(n-1)/2$.

The multiplication relations for the new bases are, for even $n$,
\begin{eqnarray}
\lefteqn{e_+^2=e_+,\; e_-^2=e_-,\; e_+e_-=0,\; e_+e_k=0,\; e_+\tilde e_k=0,\;
e_-e_k=0,\; 
e_-\tilde e_k=0,\nonumber}\\ 
&&e_k^2=e_k,\; \tilde e_k^2=-e_k,\; e_k \tilde e_k=\tilde e_k ,\; e_ke_l=0,\;
e_k\tilde e_l=0,\; \tilde e_k\tilde e_l=0,\; k\not=l,\; 
\label{e12a}
\end{eqnarray}
where $k,l=1,...,n/2-1$.
For odd $n$ the multiplication relations are
\begin{eqnarray}
\lefteqn{e_+^2=e_+,\;  e_+e_k=0,\; e_+\tilde e_k=0,\; \nonumber}\\ 
&&e_k^2=e_k,\; \tilde e_k^2=-e_k,\; e_k \tilde e_k=\tilde e_k ,\; e_ke_l=0,\;
e_k\tilde e_l=0,\; \tilde e_k\tilde e_l=0,\; k\not=l,\; 
\label{e12b}
\end{eqnarray}
where $k,l=1,...,(n-1)/2$.
The moduli of the new bases are
\begin{equation}
|e_+|=\frac{1}{\sqrt{n}},\; |e_-|=\frac{1}{\sqrt{n}},\; 
|e_k|=\sqrt{\frac{2}{n}},\; |\tilde e_k|=\sqrt{\frac{2}{n}}.
\label{e12c}
\end{equation}
It can be shown that, for even $n$,
\begin{eqnarray}
x_0+h_1x_1+\cdots+h_{n-1}x_{n-1}=e_+v_+ + e_-v_- 
+\sum_{k=1}^{n/2-1}(e_k v_k+\tilde e_k \tilde v_k),
\label{e13a}
\end{eqnarray}
and for odd $n$
\begin{eqnarray}
x_0+h_1x_1+\cdots+h_{n-1}x_{n-1}=e_+v_+ 
+\sum_{k=1}^{(n-1)/2}(e_k v_k+\tilde e_k \tilde v_k).
\label{e13b}
\end{eqnarray}
The relations (\ref{e13a}),(\ref{e13b}) give the canonical form
of a polar n-complex number.

Using the properties of the bases in Eqs. (\ref{e12a}) and (\ref{e12b}) it can
be shown that
\begin{equation}
\exp(\tilde e_k\phi_k)=1-e_k+e_k\cos\phi_k+\tilde e_k\sin\phi_k ,
\label{46a}
\end{equation}
\begin{equation}
\exp(e_k\ln\rho_k)=1-e_k+e_k\rho_k ,
\label{46b}
\end{equation}
\begin{equation}
\exp(e_+\ln v_+)=1-e_++e_+v_+
\label{46c}
\end{equation}
and, for even $n$,
\begin{equation}
\exp(e_-\ln v_-)=1-e_- +e_-v_- .
\label{46d}
\end{equation}
In Eq. (\ref{46c}), $\ln v_+$ exists as a real function 
provided that $v_+=x_0+x_1+\cdots+x_{n-1}>0$, which means that
$0<\theta_+<\pi/2$, 
and for even $n$, $\ln v_-$
exists in Eq. (\ref{46d}) as a real function provided that
$v_-=x_0-x_1+\cdots+x_{n-2}-x_{n-1}>0$, which means that $0<\theta_-<\pi/2$.
By multiplying the relations (\ref{46a})-(\ref{46d}) it results, for even $n$,
that 
\begin{equation}
\exp\left[e_+\ln v_++e_-\ln v_-+\sum_{k=1}^{n/2-1}
(e_k\ln \rho_k+\tilde e_k\phi_k)\right] 
=e_+ v_+ +e_- v_- +\sum_{k=1}^{n/2-1}
(e_k v_k+\tilde e_k \tilde v_k),
\label{47a}
\end{equation}
where the fact has ben used that 
\begin{equation}
e_+ +e_- +\sum_{k=1}^{n/2-1} e_k =1, 
\label{47b}
\end{equation}
the latter relation being a consequence of Eqs. (\ref{e11}) and (\ref{15}).
Similarly, by multiplying the relations (\ref{46a})-(\ref{46c}) it results, for
odd $n$, that
\begin{equation}
\exp\left[e_+\ln v_+ +\sum_{k=1}^{(n-1)/2}
(e_k\ln \rho_k+\tilde e_k\phi_k)\right] 
=e_+ v_+  +\sum_{k=1}^{(n-1)/2}
(e_k v_k+\tilde e_k \tilde v_k),
\label{48a}
\end{equation}
where the fact has ben used that 
\begin{equation}
e_+ +\sum_{k=1}^{(n-1)/2} e_k =1, 
\label{48b}
\end{equation}
the latter relation being a consequence of Eqs. (\ref{e12}) and (\ref{16}).

By comparing Eqs. (\ref{e13a}) and (\ref{47a}), it can be seen that, for even
$n$, 
\begin{eqnarray}
x_0+h_1x_1+\cdots+h_{n-1}x_{n-1}=
\exp\left[e_+\ln v_++e_-\ln v_-+\sum_{k=1}^{n/2-1}
(e_k\ln \rho_k+\tilde e_k\phi_k)\right] ,
\label{49a}
\end{eqnarray}
and by comparing Eqs. (\ref{e13b}) and (\ref{48a}), it can be seen that, for
odd $n$,
\begin{eqnarray}
x_0+h_1x_1+\cdots+h_{n-1}x_{n-1}
=\exp\left[e_+\ln v_+ +\sum_{k=1}^{(n-1)/2}
(e_k\ln \rho_k+\tilde e_k\phi_k)\right] .
\label{49b}
\end{eqnarray}
Using the expression of the bases in Eqs. (\ref{e11}) and (\ref{e12}) yields,
for even values of $n$, the exponential form of the n-complex number
$u=x_0+h_1x_1+\cdots+h_{n-1}x_{n-1}$ as
\begin{eqnarray}\lefteqn{
u=\rho\exp\left\{\sum_{p=1}^{n-1}h_p\left[
\frac{1}{n}\ln\frac{\sqrt{2}}{\tan\theta_+}
+\frac{(-1)^p}{n}\ln\frac{\sqrt{2}}{\tan\theta_-}\right.\right.\nonumber}\\
&&
\left.\left.-\frac{2}{n}\sum_{k=2}^{n/2-1}
\cos\left(\frac{2\pi kp}{n}\right)\ln\tan\psi_{k-1}
\right]
+\sum_{k=1}^{n/2-1}\tilde e_k\phi_k 
\right\},
\label{50a}
\end{eqnarray}
where $\rho$ is the amplitude defined in Eq. (\ref{6b}), which for even $n$ has
according to Eq. (\ref{9c}) the expression
\begin{equation}
\rho=\left(v_+v_-\rho_1^2\cdots \rho_{n/2-1}^2\right)^{1/n}.
\label{50aa}
\end{equation}
For odd values of $n$, the exponential form of the n-complex number $u$ is
\begin{eqnarray}
\lefteqn{
u=\rho\exp\left\{\sum_{p=1}^{n-1}h_p\left[
\frac{1}{n}\ln\frac{\sqrt{2}}{\tan\theta_+}
\right.\right.
\left.\left.-\frac{2}{n}\sum_{k=2}^{(n-1)/2}
\cos\left(\frac{2\pi kp}{n}\right)\ln\tan\psi_{k-1}
\right]
+\sum_{k=1}^{(n-1)/2}\tilde e_k\phi_k
\right\},\nonumber}\\
&&
\label{50b}
\end{eqnarray}
where for odd $n$, $\rho$ has according to Eq. (\ref{9d}) the expression
\begin{equation}
\rho=\left(v_+\rho_1^2\cdots \rho_{(n-1)/2}^2\right)^{1/n}.
\label{50bb}
\end{equation}

It can be checked with the aid of Eq. (\ref{46a}) that
the n-complex number $u$ can also be written, for even $n$, as
\begin{eqnarray}
x_0+h_1x_1+\cdots+h_{n-1}x_{n-1}
=\left(e_+ v_++e_- v_-+\sum_{k=1}^{n/2-1}e_k \rho_k\right)
\exp\left(\sum_{k=1}^{n/2-1}\tilde e_k\phi_k\right),
\label{51a}
\end{eqnarray}
and for odd $n$, as
\begin{eqnarray}
x_0+h_1x_1+\cdots+h_{n-1}x_{n-1}
=\left(e_+ v_++\sum_{k=1}^{(n-1)/2}e_k \rho_k\right)
\exp\left(\sum_{k=1}^{(n-1)/2}\tilde e_k\phi_k\right).
\label{51b}
\end{eqnarray}
Writing in Eqs. (\ref{51a}) and (\ref{51b}) the radius $\rho_1$, Eqs.
(\ref{20b}) and (\ref{20c}), as a factor and expressing the 
variables in terms of the polar and planar
angles with the aid of Eqs. (\ref{19b})-(\ref{19d}) yields the
trigonometric form of the n-complex number $u$, for even $n$, as
\begin{eqnarray}
\lefteqn{u=d
\left(\frac{n}{2}\right)^{1/2}
\left(\frac{1}{\tan^2\theta_+}+\frac{1}{\tan^2\theta_-}+1
+\frac{1}{\tan^2\psi_1}+\frac{1}{\tan^2\psi_2}+\cdots
+\frac{1}{\tan^2\psi_{n/2-2}}\right)^{-1/2}\nonumber}\\
&&\left(\frac{e_+\sqrt{2}}{\tan\theta_+}+\frac{e_-\sqrt{2}}{\tan\theta_-}
+e_1+\sum_{k=2}^{n/2-1}\frac{e_k}{\tan\psi_{k-1}}\right)
\exp\left(\sum_{k=1}^{n/2-1}\tilde e_k\phi_k\right),
\label{52a}
\end{eqnarray}
and for odd $n$ as
\begin{eqnarray}
\lefteqn{u=d
\left(\frac{n}{2}\right)^{1/2}
\left(\frac{1}{\tan^2\theta_+}+1
+\frac{1}{\tan^2\psi_1}+\frac{1}{\tan^2\psi_2}+\cdots
+\frac{1}{\tan^2\psi_{(n-3)/2}}\right)^{-1/2}\nonumber}\\
&&\left(\frac{e_+\sqrt{2}}{\tan\theta_+}
+e_1+\sum_{k=2}^{(n-1)/2}\frac{e_k}{\tan\psi_{k-1}}\right)
\exp\left(\sum_{k=1}^{(n-1)/2}\tilde e_k\phi_k\right).
\label{52b}
\end{eqnarray}
In Eqs. (\ref{52a}) and {\ref{52b}),  the n-complex
number $u$, written in trigonometric form, is the
product of the modulus $d$, of a part depending on the polar and planar
angles $\theta_+, \theta_-,\psi_1,...,\psi_{[(n-3)/2]},$ and of a factor
depending 
on the azimuthal angles $\phi_1,...,\phi_{[(n-1)/2]}$. 
Although the modulus of a product of n-complex numbers is not equal in
general to the product of the moduli of the factors,
it can be checked that the modulus of the factor in Eq. (\ref{52a}) is
\begin{eqnarray}
\lefteqn{
\left|\frac{e_+\sqrt{2}}{\tan\theta_+}+\frac{e_-\sqrt{2}}{\tan\theta_-}
+e_1+\sum_{k=2}^{n/2-1}\frac{e_k}{\tan\psi_{k-1}}\right|\nonumber}\\
&&=\left(\frac{2}{n}\right)^{1/2}
\left(\frac{1}{\tan^2\theta_+}+\frac{1}{\tan^2\theta_-}+1
+\frac{1}{\tan^2\psi_1}+\frac{1}{\tan^2\psi_2}+\cdots
+\frac{1}{\tan^2\psi_{n/2-2}}\right)^{1/2},\nonumber\\
&&
\label{52c}
\end{eqnarray}
and the modulus of the factor in Eq. (\ref{52b}) is
\begin{eqnarray}
\lefteqn{
\left|\frac{e_+\sqrt{2}}{\tan\theta_+}
+e_1+\sum_{k=2}^{(n-1)/2}\frac{e_k}{\tan\psi_{k-1}}\right|\nonumber}\\
&&=\left(\frac{2}{n}\right)^{1/2}
\left(\frac{1}{\tan^2\theta_+}+1
+\frac{1}{\tan^2\psi_1}+\frac{1}{\tan^2\psi_2}+\cdots
+\frac{1}{\tan^2\psi_{(n-3)/2}}\right)^{1/2}.
\label{52d}
\end{eqnarray}
Moreover, it can be checked that
\begin{eqnarray}
\left|\exp\left[\sum_{k=1}^{[(n-1)/2]}\tilde e_k\phi_k\right]\right|=1.
\label{52e}
\end{eqnarray}

The modulus $d$ in Eqs. (\ref{52a}) and (\ref{52b}) can be expressed in terms
of the amplitude $\rho$, for even $n$, as
\begin{eqnarray}
\lefteqn{d=\rho \frac{2^{(n-2)/2n}}{\sqrt{n}}
\left(\tan\theta_+\tan\theta_-
\tan^2\psi_1\cdots\tan^2\psi_{n/2-2}\right)^{1/n}\nonumber}\\
&&\left(\frac{1}{\tan^2\theta_+}+\frac{1}{\tan^2\theta_-}+1
+\frac{1}{\tan^2\psi_1}+\frac{1}{\tan^2\psi_2}+\cdots
+\frac{1}{\tan^2\psi_{n/2-2}}\right)^{1/2},
\label{53a}
\end{eqnarray}
and for odd $n$ as
\begin{eqnarray}
\lefteqn{d=\rho \frac{2^{(n-1)/2n}}{\sqrt{n}}
\left(\tan\theta_+
\tan^2\psi_1\cdots\tan^2\psi_{(n-3)/2}\right)^{1/n}\nonumber}\\
&&\left(\frac{1}{\tan^2\theta_+}+1
+\frac{1}{\tan^2\psi_1}+\frac{1}{\tan^2\psi_2}+\cdots
+\frac{1}{\tan^2\psi_{(n-3)/2}}\right)^{1/2}.
\label{53b}
\end{eqnarray}

\section{Elementary functions of a polar n-complex variable}

The logarithm $u_1$ of the n-complex number $u$, $u_1=\ln u$, can be defined
as the solution of the equation
\begin{equation}
u=e^{u_1} .
\label{54}
\end{equation}
For even $n$ the relation (\ref{47a}) shows that $\ln u$ exists 
as an n-complex function with real components if $v_+=x_0+x_1+\cdots+x_{n-1}>0$
and $v_-=x_0-x_1+\cdots+x_{n-2}-x_{n-1}>0$, which means that $0<\theta_+<\pi/2,
0<\theta_-<\pi/2$.
For odd $n$ the relation (\ref{48a}) shows that $\ln u$ exists 
as an n-complex function with real components if
$v_+=x_0+x_1+\cdots+x_{n-1}>0$, which means that $0<\theta_+<\pi/2$. 
The expression of the logarithm, obtained from Eqs. (\ref{49a}) and
(\ref{49b}), is, for even $n$,
\begin{equation}
\ln u=e_+\ln v_++e_-\ln v_-+\sum_{k=1}^{n/2-1}
(e_k\ln \rho_k+\tilde e_k\phi_k),
\label{55a}
\end{equation}
and for odd $n$ the expression is
\begin{equation}
\ln u=e_+\ln v_+ +\sum_{k=1}^{(n-1)/2}
(e_k\ln \rho_k+\tilde e_k\phi_k).
\label{55b}
\end{equation}
An expression of the logarithm depending on the amplitude $\rho$ can be
obtained from the exponential forms in Eqs. (\ref{50a}) and (\ref{50b}), for
even $n$, as
\begin{eqnarray}\lefteqn{
\ln u=\ln \rho+\sum_{p=1}^{n-1}h_p\left[
\frac{1}{n}\ln\frac{\sqrt{2}}{\tan\theta_+}
+\frac{(-1)^p}{n}\ln\frac{\sqrt{2}}{\tan\theta_-}\right.
\left.-\frac{2}{n}\sum_{k=2}^{n/2-1}
\cos\left(\frac{2\pi kp}{n}\right)\ln\tan\psi_{k-1}
\right]\nonumber}\\
&&+\sum_{k=1}^{n/2-1}\tilde e_k\phi_k,
\label{56a}
\end{eqnarray}
and for odd $n$ as
\begin{eqnarray}
\lefteqn{\ln u=\ln\rho+\sum_{p=1}^{n-1}h_p\left[
\frac{1}{n}\ln\frac{\sqrt{2}}{\tan\theta_+}
-\frac{2}{n}\sum_{k=2}^{(n-1)/2}
\cos\left(\frac{2\pi kp}{n}\right)\ln\tan\psi_{k-1}
\right]
+\sum_{k=1}^{(n-1)/2}\tilde e_k\phi_k.\nonumber}\\
&&
\label{56b}
\end{eqnarray}

The function $\ln u$ is multivalued because of the presence of the terms 
$\tilde e_k\phi_k$.
It can be inferred from Eqs. (\ref{21a})-(\ref{21g}) and (\ref{24}) that
\begin{equation}
\ln(uu^\prime)=\ln u+\ln u^\prime ,
\label{57}
\end{equation}
up to integer multiples of $2\pi\tilde e_k, k=1,...,[(n-1)/2]$.

The power function $u^m$ can be defined for real values of $m$ as
\begin{equation}
u^m=e^{m\ln u} .
\label{58}
\end{equation}
Using the expression of $\ln u$ in Eqs. (\ref{55a}) and (\ref{55b}) yields, for
even values of $n$,
\begin{equation}
u^m=e_+ v_+^m+e_- v_-^m +\sum_{k=1}^{n/2-1}
\rho_k^m(e_k\cos m\phi_k+\tilde e_k\sin m\phi_k),
\label{59a}
\end{equation}
and for odd values of $n$
\begin{equation}
u^m=e_+ v_+^m +\sum_{k=1}^{(n-1)/2}
\rho_k^m(e_k\cos m\phi_k+\tilde e_k\sin m\phi_k).
\label{59b}
\end{equation}
For integer values of $m$, the relations (\ref{59a}) and (\ref{59b}) are valid 
for any $x_0,...,x_{n-1}$.
The power function is multivalued unless $m$ is an integer. 
For integer $m$, it can be inferred from Eq. (\ref{57}) that
\begin{equation}
(uu^\prime)^m=u^m\:u^{\prime m} .
\label{59}
\end{equation}

The trigonometric functions $\cos u$ and $\sin u$ of an n-complex variable $u$
are defined by the series
\begin{equation}
\cos u = 1 - u^2/2!+u^4/4!+\cdots, 
\label{60}
\end{equation}
\begin{equation}
\sin u=u-u^3/3!+u^5/5! +\cdots .
\label{61}
\end{equation}
It can be checked by series multiplication that the usual addition theorems
hold for the n-complex numbers $u, u^\prime$,
\begin{equation}
\cos(u+u^\prime)=\cos u\cos u^\prime - \sin u \sin u^\prime ,
\label{62}
\end{equation}
\begin{equation}
\sin(u+u^\prime)=\sin u\cos u^\prime + \cos u \sin u^\prime .
\label{63}
\end{equation}

In order to obtain expressions for the trigonometric functions of n-complex
variables, these will be expressed with the aid of the imaginary unit $i$ as
\begin{equation}
\cos u=\frac{1}{2}(e^{iu}+e^{-iu}),\:\sin u=\frac{1}{2i}(e^{iu}-e^{-iu}).
\label{64}
\end{equation}
The imaginary unit $i$ is used for the convenience of notations, and it does
not appear in the final results.
The validity of Eq. (\ref{64}) can be checked by comparing the series for the
two sides of the relations.
Since the expression of the exponential function $e^{h_k y}$ in terms of the
units $1, h_1, ... h_{n-1}$ given in Eq. (\ref{28b}) depends on the polar
cosexponential functions $g_{np}(y)$, the expression of the trigonometric
functions will depend on the functions
$g_{p+}^{(c)}(y)=(1/2)[g_{np}(iy)+g_{np}(-iy)]$ and 
$g_{p-}^{(c)}(y)=(1/2i)[g_{np}(iy)-g_{np}(-iy)]$, 
\begin{equation}
\cos(h_k y)=\sum_{p=0}^{n-1}h_{kp-n[kp/n]}g_{p+}^{(c)}(y),
\label{66a}
\end{equation}
\begin{equation}
\sin(h_k y)=\sum_{p=0}^{n-1}h_{kp-n[kp/n]}g_{p-}^{(c)}(y),
\label{66b}
\end{equation}
where
\begin{eqnarray}
\lefteqn{g_{p+}^{(c)}(y)=\frac{1}{n}\sum_{l=0}^{n-1}\left\{
\cos\left[y\cos\left(\frac{2\pi l}{n}\right)\right]
\cosh\left[y\sin\left(\frac{2\pi l}{n}\right)\right]
\cos\left(\frac{2\pi lp}{n}\right)\right.\nonumber}\\
&&\left.-\sin\left[y\cos\left(\frac{2\pi l}{n}\right)\right]
\sinh\left[y\sin\left(\frac{2\pi l}{n}\right)\right]
\sin\left(\frac{2\pi lp}{n}\right)
\right\},
\label{65a}
\end{eqnarray}
\begin{eqnarray}
\lefteqn{g_{p-}^{(c)}(y)=\frac{1}{n}\sum_{l=0}^{n-1}\left\{
\sin\left[y\cos\left(\frac{2\pi l}{n}\right)\right]
\cosh\left[y\sin\left(\frac{2\pi l}{n}\right)\right]
\cos\left(\frac{2\pi lp}{n}\right)\right.\nonumber}\\
&&\left.+\cos\left[y\cos\left(\frac{2\pi l}{n}\right)\right]
\sinh\left[y\sin\left(\frac{2\pi l}{n}\right)\right]
\sin\left(\frac{2\pi lp}{n}\right)
\right\}.
\label{65b}
\end{eqnarray}

The hyperbolic functions $\cosh u$ and $\sinh u $ of the n-complex variable
$u$ can be defined by the series
\begin{equation}
\cosh u = 1 + u^2/2!+u^4/4!+\cdots, 
\label{66}
\end{equation}
\begin{equation}
\sinh u=u+u^3/3!+u^5/5! +\cdots .
\label{67}
\end{equation}
It can be checked by series multiplication that the usual addition theorems
hold for the n-complex numbers $u, u^\prime$,
\begin{equation}
\cosh(u+u^\prime)=\cosh u\cosh u^\prime + \sinh u \sinh u^\prime ,
\label{68}
\end{equation}
\begin{equation}
\sinh(u+u^\prime)=\sinh u\cosh u^\prime + \cosh u \sinh u^\prime .
\label{69}
\end{equation}
In order to obtain expressions for the hyperbolic functions of n-complex
variables, these will be expressed as
\begin{equation}
\cosh u=\frac{1}{2}(e^{u}+e^{-u}),\:\sinh u=\frac{1}{2}(e^{u}-e^{-u}).
\label{70}
\end{equation}
The validity of Eq. (\ref{70}) can be checked by comparing the series for the
two sides of the relations.
Since the expression of the exponential function $e^{h_k y}$ in terms of the
units $1, h_1, ... h_{n-1}$ given in Eq. (\ref{28b}) depends on the polar
cosexponential functions $g_{np}(y)$, the expression of the hyperbolic
functions will depend on the even part $g_{p+}(y)=(1/2)[g_{np}(y)+g_{np}(-y)]$
and on 
the odd part $g_{p-}(y)=(1/2)[g_{np}(y)-g_{np}(-y)]$ of $g_{np}$, 
\begin{equation}
\cosh(h_k y)=\sum_{p=0}^{n-1}h_{kp-n[kp/n]}g_{p+}(y),
\label{71a}
\end{equation}
\begin{equation}
\sinh(h_k y)=\sum_{p=0}^{n-1}h_{kp-n[kp/n]}g_{p-}(y),
\label{71b}
\end{equation}
where
\begin{eqnarray}
\lefteqn{g_{p+}(y)=\frac{1}{n}\sum_{l=0}^{n-1}\left\{
\cosh\left[y\cos\left(\frac{2\pi l}{n}\right)\right]
\cos\left[y\sin\left(\frac{2\pi l}{n}\right)\right]
\cos\left(\frac{2\pi lp}{n}\right)\right.\nonumber}\\
&&\left.+\sinh\left[y\cos\left(\frac{2\pi l}{n}\right)\right]
\sin\left[y\sin\left(\frac{2\pi l}{n}\right)\right]
\sin\left(\frac{2\pi lp}{n}\right)
\right\},
\label{72a}
\end{eqnarray}
\begin{eqnarray}
\lefteqn{g_{p-}(y)=\frac{1}{n}\sum_{l=0}^{n-1}\left\{
\sinh\left[y\cos\left(\frac{2\pi l}{n}\right)\right]
\cos\left[y\sin\left(\frac{2\pi l}{n}\right)\right]
\cos\left(\frac{2\pi lp}{n}\right)\right.\nonumber}\\
&&\left.+\cosh\left[y\cos\left(\frac{2\pi l}{n}\right)\right]
\sin\left[y\sin\left(\frac{2\pi l}{n}\right)\right]
\sin\left(\frac{2\pi lp}{n}\right)
\right\}.
\label{72b}
\end{eqnarray}

The exponential, trigonometric and hyperbolic functions can also be expressed
with the aid of the bases introduced in Eqs. (\ref{e11}) and (\ref{e12}).
Using the expression of the n-complex number in Eq. (\ref{e13a}), for even $n$,
yields for the exponential of the n-complex variable $u$
\begin{eqnarray}
e^u=e_+e^{v_+} + e_-e^{v_-} 
+\sum_{k=1}^{n/2-1}e^{v_k}\left(e_k \cos \tilde v_k+\tilde e_k \sin\tilde
v_k\right).
\label{73a}
\end{eqnarray}
For odd $n$, the expression of the n-complex variable in Eq. (\ref{e13b})
yileds for the exponential 
\begin{eqnarray}
e^u=e_+e^{v_+}  
+\sum_{k=1}^{(n-1)/2}e^{v_k}\left(e_k \cos \tilde v_k+\tilde e_k \sin\tilde
v_k\right).
\label{73b}
\end{eqnarray}

The trigonometric functions can be obtained from Eqs. (\ref{73a}) and
(\ref{73b} with the aid of Eqs. (\ref{64}). The trigonometric functions of the
n-complex variable $u$ are, for even $n$,
\begin{equation}
\cos u=e_+\cos v_+ + e_-\cos v_- 
+\sum_{k=1}^{n/2-1}\left(e_k \cos v_k\cosh \tilde v_k
-\tilde e_k \sin v_k\sinh\tilde v_k\right),
\label{74a}
\end{equation}
\begin{equation}
\sin u=e_+\sin v_+ + e_-\sin v_- 
+\sum_{k=1}^{n/2-1}\left(e_k \sin v_k\cosh \tilde v_k
+\tilde e_k \cos v_k\sinh\tilde v_k\right),
\label{74b}
\end{equation}
and for odd $n$ the trigonometric functions are
\begin{equation}
\cos u=e_+\cos v_+  
+\sum_{k=1}^{(n-1)/2}\left(e_k \cos v_k\cosh \tilde v_k
-\tilde e_k \sin v_k\sinh\tilde v_k\right),
\label{74c}
\end{equation}
\begin{equation}
\sin u=e_+\sin v_+  
+\sum_{k=1}^{(n-1)/2}\left(e_k \sin v_k\cosh \tilde v_k
+\tilde e_k \cos v_k\sinh\tilde v_k\right).
\label{74d}
\end{equation}

The hyperbolic functions can be obtained from Eqs. (\ref{73a}) and
(\ref{73b} with the aid of Eqs. (\ref{70}). The hyperbolic functions of the
n-complex variable $u$ are, for even $n$,
\begin{equation}
\cosh u=e_+\cosh v_+ + e_-\cosh v_- 
+\sum_{k=1}^{n/2-1}\left(e_k \cosh v_k\cos \tilde v_k
+\tilde e_k \sinh v_k\sin\tilde v_k\right),
\label{75a}
\end{equation}
\begin{equation}
\sinh u=e_+\sinh v_+ + e_-\sinh v_- 
+\sum_{k=1}^{n/2-1}\left(e_k \sinh v_k\cos \tilde v_k
+\tilde e_k \cosh v_k\sin\tilde v_k\right),
\label{75b}
\end{equation}
and for odd $n$ the hyperbolic functions are
\begin{equation}
\cosh u=e_+\cosh v_+  
+\sum_{k=1}^{(n-1)/2}\left(e_k \cosh v_k\cos \tilde v_k
+\tilde e_k \sinh v_k\sin\tilde v_k\right),
\label{75c}
\end{equation}
\begin{equation}
\sinh u=e_+\sinh v_+  
+\sum_{k=1}^{(n-1)/2}\left(e_k \sinh v_k\cos \tilde v_k
+\tilde e_k \cosh v_k\sin\tilde v_k\right).
\label{75d}
\end{equation}

\section{Power series of polar n-complex numbers}

An n-complex series is an infinite sum of the form
\begin{equation}
a_0+a_1+a_2+\cdots+a_n+\cdots , 
\label{76}
\end{equation}
where the coefficients $a_n$ are n-complex numbers. The convergence of 
the series (\ref{76}) can be defined in terms of the convergence of its $n$
real components. The convergence of a n-complex series can also be studied
using n-complex variables. The main criterion for absolute convergence 
remains the comparison theorem, but this requires a number of inequalities
which will be discussed further.

The modulus $d=|u|$ of an n-complex number $u$ has been defined in Eq.
(\ref{10}). Since $|x_0|\leq |u|, |x_1|\leq |u|,..., |x_{n-1}|\leq |u|$, a
property of absolute convergence established via a comparison theorem based on
the modulus of the series (\ref{76}) will ensure the absolute convergence of
each real component of that series.

The modulus of the sum $u_1+u_2$ of the n-complex numbers $u_1, u_2$ fulfils
the inequality
\begin{equation}
||u^\prime|-|u^{\prime\prime}||\leq |u^\prime+u^{\prime\prime}|\leq
|u^\prime|+|u^{\prime\prime}| . 
\label{78}
\end{equation}
For the product, the relation is 
\begin{equation}
|u^\prime u^{\prime\prime}|\leq \sqrt{n}|u^\prime||u^{\prime\prime}| ,
\label{79}
\end{equation}
as can be shown from Eqs. (\ref{17}) and (\ref{18}). The relation (\ref{79})
replaces the relation of equality extant between 2-dimensional regular complex
numbers. The equality in Eq. (\ref{79}) takes place for
$\rho_1\rho_1^\prime=0,..., \rho_{[(n-1)/2]}\rho_{[(n-1)/2]}^\prime=0$ and, for
even $n$, for $v_+v_-^\prime=0$, $v_-v_+^\prime=0$.  

From Eq. (\ref{79}) it results, for $u=u^\prime$, that
\begin{equation}
|u^2|\leq \sqrt{n} |u|^2 .
\label{80}
\end{equation}
The relation in Eq. (\ref{80}) becomes an equality for 
$\rho_1=0,...,\rho_{[(n-1)/2]}=0$ and, for even $n$, $v_+=0$ or $v_-=0$.
The inequality in Eq. (\ref{79}) implies that
\begin{equation}
|u^l|\leq n^{(l-1)/2}|u|^l ,
\label{81}
\end{equation}
where $l$ is a natural number.
From Eqs. (\ref{79}) and (\ref{81}) it results that
\begin{equation}
|au^l|\leq n^{l/2} |a| |u|^l .
\label{82}
\end{equation}

A power series of the n-complex variable $u$ is a series of the form
\begin{equation}
a_0+a_1 u + a_2 u^2+\cdots +a_l u^l+\cdots .
\label{83}
\end{equation}
Since
\begin{equation}
\left|\sum_{l=0}^\infty a_l u^l\right| \leq  \sum_{l=0}^\infty
n^{l/2} |a_l| |u|^l ,
\label{84}
\end{equation}
a sufficient condition for the absolute convergence of this series is that
\begin{equation}
\lim_{l\rightarrow \infty}\frac{\sqrt{n}|a_{l+1}||u|}{|a_l|}<1 .
\label{85}
\end{equation}
Thus the series is absolutely convergent for 
\begin{equation}
|u|<c,
\label{86}
\end{equation}
where 
\begin{equation}
c=\lim_{l\rightarrow\infty} \frac{|a_l|}{\sqrt{n}|a_{l+1}|} .
\label{87}
\end{equation}

The convergence of the series (\ref{83}) can be also studied with the aid of
the formulas (\ref{59a}), (\ref{59b}) which for integer values of $m$ are
valid for any values of $x_0,...,x_{n-1}$, as mentioned previously.
If $a_l=\sum_{p=0}^{n-1}h_p a_{lp}$, and
\begin{equation}
A_{l+}=\sum_{p=0}^{n-1}a_{lp},
\label{88a}
\end{equation}
\begin{equation}
A_{lk}=\sum_{p=0}^{n-1}a_{lp}\cos\frac{2\pi kp}{n},
\label{88b}
\end{equation}
\begin{equation}
\tilde A_{lk}=\sum_{p=0}^{n-1}a_{lp}\sin\frac{2\pi kp}{n},
\label{88c}
\end{equation}
for $k=1,...,[(n-1)/2]$, and for even $n$ 
\begin{equation}
A_{l-}=\sum_{p=0}^{n-1}(-1)^p a_{lp},
\label{88d}
\end{equation}
the series (\ref{83}) can be written, for even $n$, as
\begin{equation}
\sum_{l=0}^\infty \left[
e_+A_{l+}v_+^l+e_-A_{l-}v_-^l+\sum_{k=1}^{n/2-1}
(e_k A_{lk}+\tilde e_k\tilde A_{lk})(e_k v_k+\tilde e_k\tilde v_k)^l 
\right],
\label{89a}
\end{equation}
and for odd $n$ as
\begin{equation}
\sum_{l=0}^\infty \left[
e_+A_{l+}v_+^l+\sum_{k=1}^{(n-1)/2}
(e_k A_{lk}+\tilde e_k\tilde A_{lk})(e_k v_k+\tilde e_k\tilde v_k)^l 
\right].
\label{89b}
\end{equation}

The series in Eq. (\ref{83}) is absolutely convergent for 
\begin{equation}
|v_+|<c_+,\:
|v_-|<c_-,\:
\rho_k<c_k, 
\label{90}
\end{equation}
for $k=1,..., [(n-1)/2]$, where 
\begin{equation}
c_+=\lim_{l\rightarrow\infty} \frac{|A_{l+}|}{|A_{l+1,+}|} ,\:
c_-=\lim_{l\rightarrow\infty} \frac{|A_{l-}|}{|A_{l+1,-}|} ,\:
c_k=\lim_{l\rightarrow\infty} \frac
{\left(A_{lk}^2+\tilde A_{lk}^2\right)^{1/2}}
{\left(A_{l+1,k}^2+\tilde A_{l+1,k}^2\right)^{1/2}} .
\label{91}
\end{equation}
The relations (\ref{90}) show that the region of convergence of the series
(\ref{83}) is an n-dimensional cylinder.

It can be shown that, for even $n$, $c=(1/\sqrt{n})\;{\rm
min}(c_+,c_-,c_1,...,c_{n/2-1})$, and for odd $n$ $c=(1/\sqrt{n})\;{\rm
min}(c_+,c_1,...,c_{(n-1)/2})$, where ${\rm min}$ designates the smallest of
the numbers in the argument of this function. Using the expression of $|u|$ in
Eqs. 
(\ref{17}) or (\ref{18}), it can be seen that the spherical region of
convergence defined in Eqs. (\ref{86}), (\ref{87}) is a subset of the
cylindrical region of convergence defined in Eqs. (\ref{90}) and (\ref{91}).

\section{Analytic functions of polar n-complex variables}

The derivative  
of a function $f(u)$ of the n-complex variables $u$ is
defined as a function $f^\prime (u)$ having the property that
\begin{equation}
|f(u)-f(u_0)-f^\prime (u_0)(u-u_0)|\rightarrow 0 \:\:{\rm as} 
\:\:|u-u_0|\rightarrow 0 . 
\label{h88}
\end{equation}
If the difference $u-u_0$ is not parallel to one of the nodal hypersurfaces,
the definition in Eq. (\ref{h88}) can also 
be written as
\begin{equation}
f^\prime (u_0)=\lim_{u\rightarrow u_0}\frac{f(u)-f(u_0)}{u-u_0} .
\label{h89}
\end{equation}
The derivative of the function $f(u)=u^m $, with $m$ an integer, 
is $f^\prime (u)=mu^{m-1}$, as can be seen by developing $u^m=[u_0+(u-u_0)]^m$
as
\begin{equation}
u^m=\sum_{p=0}^{m}\frac{m!}{p!(m-p)!}u_0^{m-p}(u-u_0)^p,
\label{h90}
\end{equation}
and using the definition (\ref{h88}).

If the function $f^\prime (u)$ defined in Eq. (\ref{h88}) is independent of the
direction in space along which $u$ is approaching $u_0$, the function $f(u)$ 
is said to be analytic, analogously to the case of functions of regular complex
variables. \cite{3} 
The function $u^m$, with $m$ an integer, 
of the n-complex variable $u$ is analytic, because the
difference $u^m-u_0^m$ is always proportional to $u-u_0$, as can be seen from
Eq. (\ref{h90}). Then series of
integer powers of $u$ will also be analytic functions of the n-complex
variable $u$, and this result holds in fact for any commutative algebra. 

If an analytic function is defined by a series around a certain point, for
example $u=0$, as
\begin{equation}
f(u)=\sum_{k=0}^\infty a_k u^k ,
\label{h91a}
\end{equation}
an expansion of $f(u)$ around a different point $u_0$,
\begin{equation}
f(u)=\sum_{k=0}^\infty c_k (u-u_0)^k ,
\label{h91aa}
\end{equation}
can be obtained by
substituting in Eq. (\ref{h91a}) the expression of $u^k$ according to Eq.
(\ref{h90}). Assuming that the series are absolutely convergent so that the
order of the terms can be modified and ordering the terms in the resulting
expression according to the increasing powers of $u-u_0$ yields
\begin{equation}
f(u)=\sum_{k,l=0}^\infty \frac{(k+l)!}{k!l!}a_{k+l} u_0^l (u-u_0)^k .
\label{h91b}
\end{equation}
Since the derivative of order $k$ at $u=u_0$ of the function $f(u)$ , Eq.
(\ref{h91a}), is 
\begin{equation}
f^{(k)}(u_0)=\sum_{l=0}^\infty \frac{(k+l)!}{l!}a_{k+l} u_0^l ,
\label{h91c}
\end{equation}
the expansion of $f(u)$ around $u=u_0$, Eq. (\ref{h91b}), becomes
\begin{equation}
f(u)=\sum_{k=0}^\infty \frac{1}{k!} f^{(k)}(u_0)(u-u_0)^k ,
\label{h91d}
\end{equation}
which has the same form as the series expansion of 2-dimensional complex
functions. 
The relation (\ref{h91d}) shows that the coefficients in the series expansion,
Eq. (\ref{h91aa}), are
\begin{equation}
c_k=\frac{1}{k!}f^{(k)}(u_0) .
\label{h92}
\end{equation}

The rules for obtaining the derivatives and the integrals of the basic
functions can 
be obtained from the series of definitions and, as long as these series
expansions have the same form as the corresponding series for the
2-dimensional complex functions, the rules of derivation and integration remain
unchanged. 

If the n-complex function $f(u)$
of the n-complex variable $u$ is written in terms of 
the real functions $P_k(x_0,...,x_{n-1}), k=0,1,...,n-1$ of the real
variables $x_0,x_1,...,x_{n-1}$ as 
\begin{equation}
f(u)=\sum_{k=0}^{n-1}h_kP_k(x_0,...,x_{n-1}),
\label{h93}
\end{equation}
then relations of equality 
exist between the partial derivatives of the functions $P_k$. 
The derivative of the function $f$ can be written as
\begin{eqnarray}
\lim_{\Delta u\rightarrow 0}\frac{1}{\Delta u} 
\sum_{k=0}^{n-1}\left(h_k\sum_{l=0}^{n-1}
\frac{\partial P_k}{\partial x_l}\Delta x_l\right),
\label{h94}
\end{eqnarray}
where
\begin{equation}
\Delta u=\sum_{k=0}^{n-1}h_l\Delta x_l.
\label{h94a}
\end{equation}
The
relations between the partials derivatives of the functions $P_k$ are
obtained by setting successively in   
Eq. (\ref{h94}) $\Delta u=h_l\Delta x_l$, for $l=0,1,...,n-1$, and equating the
resulting expressions. 
The relations are 
\begin{equation}
\frac{\partial P_k}{\partial x_0} = \frac{\partial P_{k+1}}{\partial x_1} 
=\cdots=\frac{\partial P_{n-1}}{\partial x_{n-k-1}} 
= \frac{\partial P_0}{\partial x_{n-k}}=\cdots
=\frac{\partial P_{k-1}}{\partial x_{n-1}}, 
\label{h95}
\end{equation}
for $k=0,1,...,n-1$.
The relations (\ref{h95}) are analogous to the Riemann relations
for the real and imaginary components of a complex function. 
It can be shown from Eqs. (\ref{h95}) that the components $P_k$ fulfil the
second-order equations
\begin{eqnarray}
\lefteqn{\frac{\partial^2 P_k}{\partial x_0\partial x_l}
=\frac{\partial^2 P_k}{\partial x_1\partial x_{l-1}}
=\cdots=
\frac{\partial^2 P_k}{\partial x_{[l/2]}\partial x_{l-[l/2]}}}\nonumber\\
&&=\frac{\partial^2 P_k}{\partial x_{l+1}\partial x_{n-1}}
=\frac{\partial^2 P_k}{\partial x_{l+2}\partial x_{n-2}}
=\cdots
=\frac{\partial^2 P_k}{\partial x_{l+1+[(n-l-2)/2]}
\partial x_{n-1-[(n-l-2)/2]}} ,
\label{96}
\end{eqnarray}
for $k,l=0,1,...,n-1$.

\section{Integrals of polar n-complex functions}

The singularities of n-complex functions arise from terms of the form
$1/(u-u_0)^n$, with $n>0$. Functions containing such terms are singular not
only at $u=u_0$, but also at all points of the hypersurfaces
passing through the pole $u_0$ and which are parallel to the nodal
hypersurfaces.  

The integral of an n-complex function between two points $A, B$ along a path
situated in a region free of singularities is independent of path, which means
that the integral of an analytic function along a loop situated in a region
free of singularities is zero,
\begin{equation}
\oint_\Gamma f(u) du = 0,
\label{111}
\end{equation}
where it is supposed that a surface $\Sigma$ spanning 
the closed loop $\Gamma$ is not intersected by any of
the hypersurfaces associated with the
singularities of the function $f(u)$. Using the expression, Eq. (\ref{h93}),
for $f(u)$ and the fact that 
\begin{eqnarray}
du=\sum_{k=0}^{n-1}h_k dx_k, 
\label{111a}
\end{eqnarray}
the explicit form of the integral in Eq. (\ref{111}) is
\begin{eqnarray}
\oint _\Gamma f(u) du = \oint_\Gamma
\sum_{k=0}^{n-1}h_k\sum_{l=0}^{n-1}P_l dx_{k-l+n[(n-k-1+l)/n]}.
\label{112}
\end{eqnarray}

If the functions $P_k$ are regular on a surface $\Sigma$
spanning the loop $\Gamma$,
the integral along the loop $\Gamma$ can be transformed in an integral over the
surface $\Sigma$ of terms of the form
$\partial P_l/\partial x_{k-m+n[(n-k+m-1)/n]} 
-  \partial P_m/\partial x_{k-l+n[(n-k+l-1)/n]}$.
These terms are equal to zero by Eqs. (\ref{h95}), and this
proves Eq. (\ref{111}). 

The integral of the function $(u-u_0)^m$ on a closed loop $\Gamma$ is equal to
zero for $m$ a positive or negative integer not equal to -1,
\begin{equation}
\oint_\Gamma (u-u_0)^m du = 0, \:\: m \:\:{\rm integer},\: m\not=-1 .
\label{112b}
\end{equation}
This is due to the fact that $\int (u-u_0)^m du=(u-u_0)^{m+1}/(m+1), $ and to
the fact that the function $(u-u_0)^{m+1}$ is singlevalued for $m$ an integer.

The integral $\oint_\Gamma du/(u-u_0)$ can be calculated using the exponential
form, Eqs. (\ref{50a}) and (\ref{50b}), for the difference $u-u_0$, which for
even $n$ is 
\begin{eqnarray}\lefteqn{
u-u_0=\rho\exp\left\{\sum_{p=1}^{n-1}h_p\left[
\frac{1}{n}\ln\frac{\sqrt{2}}{\tan\theta_+}
+\frac{(-1)^p}{n}\ln\frac{\sqrt{2}}{\tan\theta_-}\right.\right.\nonumber}\\
&&
\left.\left.-\frac{2}{n}\sum_{k=2}^{n/2-1}
\cos\left(\frac{2\pi kp}{n}\right)\ln\tan\psi_{k-1}
\right]
+\sum_{k=1}^{n/2-1}\tilde e_k\phi_k\right\},
\label{113a}
\end{eqnarray}
and for odd $n$ is
\begin{eqnarray}
\lefteqn{
u-u_0=\rho\exp\left\{\sum_{p=1}^{n-1}h_p\left[
\frac{1}{n}\ln\frac{\sqrt{2}}{\tan\theta_+}
\right.\right.
\left.\left.-\frac{2}{n}\sum_{k=2}^{(n-1)/2}
\cos\left(\frac{2\pi kp}{n}\right)\ln\tan\psi_{k-1}
\right]
+\sum_{k=1}^{(n-1)/2}\tilde e_k\phi_k\right\}.\nonumber}\\
&&
\label{113b}
\end{eqnarray}
Thus for even $n$ the quantity $du/(u-u_0)$ is
\begin{eqnarray}\lefteqn{
\frac{du}{u-u_0}=
\frac{d\rho}{\rho}
+\sum_{p=1}^{n-1}h_p\left[
\frac{1}{n}d\ln\frac{\sqrt{2}}{\tan\theta_+}
+\frac{(-1)^p}{n}d\ln\frac{\sqrt{2}}{\tan\theta_-}\right.\nonumber}\\
&&
\left.-\frac{2}{n}\sum_{k=2}^{n/2-1}
\cos\left(\frac{2\pi kp}{n}\right)d\ln\tan\psi_{k-1}
\right]
+\sum_{k=1}^{n/2-1}\tilde e_kd\phi_k,
\label{114a}
\end{eqnarray}
and for odd $n$
\begin{eqnarray}
\lefteqn{
\frac{du}{u-u_0}=\frac{d\rho}{\rho}
+\sum_{p=1}^{n-1}h_p\left[
\frac{1}{n}d\ln\frac{\sqrt{2}}{\tan\theta_+}
\right.
\left.-\frac{2}{n}\sum_{k=2}^{(n-1)/2}
\cos\left(\frac{2\pi kp}{n}\right)d\ln\tan\psi_{k-1}
\right]
+\sum_{k=1}^{(n-1)/2}\tilde e_kd\phi_k
.\nonumber}\\
&&
\label{114b}
\end{eqnarray}
Since $\rho, \ln(\sqrt{2}/\tan\theta_+),\ln(\sqrt{2}/\tan\theta_-),
\ln(\tan\psi_{k-1})$ are singlevalued variables, it follows that
$\oint_\Gamma d\rho/\rho =0, 
\oint_\Gamma d(\ln\sqrt{2}/\tan\theta_+)=0,
\oint_\Gamma d(\ln\sqrt{2}/\tan\theta_-)=0,
\oint_\Gamma d(\ln\tan\psi_{k-1})=0$.
On the other hand since,
$\phi_k$ are cyclic variables, they may give contributions to
the integral around the closed loop $\Gamma$.

The expression of $\oint_\Gamma du/(u-u_0)$ can be written 
with the aid of a functional which will be called int($M,C$), defined for a
point $M$ and a closed curve $C$ in a two-dimensional plane, such that 
\begin{equation}
{\rm int}(M,C)=\left\{
\begin{array}{l}
1 \;\:{\rm if} \;\:M \;\:{\rm is \;\:an \;\:interior \;\:point \;\:of} \;\:C
,\\  
0 \;\:{\rm if} \;\:M \;\:{\rm is \;\:exterior \;\:to}\:\; C .\\
\end{array}\right.
\label{118}
\end{equation}
With this notation the result of the integration on a closed path $\Gamma$
can be written as 
\begin{equation}
\oint_\Gamma\frac{du}{u-u_0}=
\sum_{k=1}^{[(n-1)/2]}2\pi\tilde e_k 
\;{\rm int}(u_{0\xi_k\eta_k},\Gamma_{\xi_k\eta_k}) ,
\label{119}
\end{equation}
where $u_{0\xi_k\eta_k}$ and $\Gamma_{\xi_k\eta_k}$ are respectively the
projections of the point $u_0$ and of 
the loop $\Gamma$ on the plane defined by the axes $\xi_k$ and $\eta_k$,
as shown in Fig. 3.

If $f(u)$ is an analytic n-complex function which can be expanded in a
series as written in Eq. (\ref{h91aa}), and the expansion holds on the curve
$\Gamma$ and on a surface spanning $\Gamma$, then from Eqs. (\ref{112b}) and
(\ref{119}) it follows that
\begin{equation}
\oint_\Gamma \frac{f(u)du}{u-u_0}=
2\pi f(u_0)\sum_{k=1}^{[(n-1)/2]}\tilde e_k 
\;{\rm int}(u_{0\xi_k\eta_k},\Gamma_{\xi_k\eta_k}) .
\label{120}
\end{equation}

Substituting in the right-hand side of 
Eq. (\ref{120}) the expression of $f(u)$ in terms of the real 
components $P_k$, Eq. (\ref{h93}), yields
\begin{equation}
\oint_\Gamma \frac{f(u)du}{u-u_0}=
\frac{2}{n}\sum_{k=1}^{[(n-1)/2]}\sum_{l,m=0}^{n-1}
h_l\sin\left[\frac{2\pi(l-m)k}{n}\right] P_m(u_0)
\;{\rm int}(u_{0\xi_k\eta_k},\Gamma_{\xi_k\eta_k}) .
\label{121}
\end{equation}
It the integral in Eq. (\ref{121}) is written as 
\begin{equation}
\oint_\Gamma \frac{f(u)du}{u-u_0}=\sum_{l=0}^{n-1}h_l I_l,
\label{122a}
\end{equation}
it can be checked that
\begin{equation}
\sum_{l=0}^{n-1} I_l=0.
\label{122b}
\end{equation}

If $f(u)$ can be expanded as written in Eq. (\ref{h91aa}) on 
$\Gamma$ and on a surface spanning $\Gamma$, then from Eqs. (\ref{112b}) and
(\ref{119}) it also results that
\begin{equation}
\oint_\Gamma \frac{f(u)du}{(u-u_0)^{n+1}}=
\frac{2\pi}{n!}f^{(n)}(u_0)\sum_{k=1}^{[(n-1)/2]}\tilde e_k 
\;{\rm int}(u_{0\xi_k\eta_k},\Gamma_{\xi_k\eta_k}) ,
\label{122}
\end{equation}
where the fact has been used  that the derivative $f^{(n)}(u_0)$ is related to
the expansion coefficient in Eq. (\ref{h91aa}) according to Eq. (\ref{h92}).

If a function $f(u)$ is expanded in positive and negative powers of $u-u_l$,
where $u_l$ are n-complex constants, $l$ being an index, the integral of $f$
on a closed loop $\Gamma$ is determined by the terms in the expansion of $f$
which are of the form $r_l/(u-u_l)$,
\begin{equation}
f(u)=\cdots+\sum_l\frac{r_l}{u-u_l}+\cdots.
\label{123}
\end{equation}
Then the integral of $f$ on a closed loop $\Gamma$ is
\begin{equation}
\oint_\Gamma f(u) du = 
2\pi \sum_l\sum_{k=1}^{[(n-1)/2]}\tilde e_k 
\;{\rm int}(u_{l\xi_k\eta_k},\Gamma_{\xi_k\eta_k})r_l .
\label{124}
\end{equation}

\section{Factorization of polar n-complex polynomials}

A polynomial of degree $m$ of the n-complex variable $u$ has the form
\begin{equation}
P_m(u)=u^m+a_1 u^{m-1}+\cdots+a_{m-1} u +a_m ,
\label{125}
\end{equation}
where $a_l$, for $l=1,...,m$, are in general n-complex constants.
If $a_l=\sum_{p=0}^{n-1}h_p a_{lp}$, and with the
notations of Eqs. (\ref{88a})-(\ref{88d}) applied for $l= 1, \cdots, m$, the
polynomial $P_m(u)$ can be written, for even $n$, as 
\begin{eqnarray}
\lefteqn{P_m= 
e_+\left(v_+^m +\sum_{l=1}^{m}A_{l+}v_+^{m-l} \right)
+e_-\left(v_-^m +\sum_{l=1}^{m}A_{l-}v_-^{m-l} \right) \nonumber}\\
&&+\sum_{k=1}^{n/2-1}
\left[(e_k v_k+\tilde e_k\tilde v_k)^m+
\sum_{l=1}^m(e_k A_{lk}+\tilde e_k\tilde A_{lk})
(e_k v_k+\tilde e_k\tilde v_k)^{m-l} 
\right],
\label{126a}
\end{eqnarray}
where the constants $A_{l+}, A_{l-}, A_{lk}, \tilde A_{lk}$ 
are real numbers.
For odd $n$ the expression of the polynomial is
\begin{eqnarray}
\lefteqn{P_m= 
e_+\left(v_+^m +\sum_{l=1}^{m}A_{l+}v_+^{m-l} \right) \nonumber}\\
&&+\sum_{k=1}^{(n-1)/2}
\left[(e_k v_k+\tilde e_k\tilde v_k)^m+
\sum_{l=1}^m(e_k A_{lk}+\tilde e_k\tilde A_{lk})
(e_k v_k+\tilde e_k\tilde v_k)^{m-l} 
\right].
\label{126b}
\end{eqnarray}

The polynomials of degree $m$ in $e_k v_k+\tilde e_k\tilde v_k$ 
in Eqs. (\ref{126a}) and (\ref{126b})
can always be written as a product of linear factors of the form
$e_k (v_k-v_{kp})+\tilde e_k(\tilde v_k-\tilde v_{kp})$, where the
constants $v_{kp}, \tilde v_{kp}$ are real.
The polynomials of degree $m$ with real coefficients in Eqs. (\ref{126a}) 
and (\ref{126b})
which are multiplied by $e_+$ and $e_-$ can be written as a product
of linear or quadratic factors with real coefficients, or as a product of
linear factors which, if imaginary, appear always in complex conjugate pairs.
Using the latter form for the simplicity of notations, the polynomial $P_m$
can be written, for even $n$, as
\begin{equation}
P_m=e_+\prod_{p=1}^m (v_+ -v_{p+})+e_-\prod_{p=1}^m (v_- -v_{p-})
+\sum_{k=1}^{n/2-1}\prod_{p=1}^m
\left\{e_k (v_k-v_{kp})+\tilde e_k(\tilde v_k-\tilde v_{kp})\right\},
\label{127a}
\end{equation}
where the quantities $v_{p+}$ appear always in complex conjugate
pairs, and the quantities $\tilde v_{p-}$ appear always in complex
conjugate pairs.
For odd $n$ the polynomial can be written as
\begin{equation}
P_m=e_+\prod_{p=1}^m (v_+ -v_{p+})
+\sum_{k=1}^{(n-1)/2}\prod_{p=1}^m
\left\{e_k (v_k-v_{kp})+\tilde e_k(\tilde v_k-\tilde v_{kp})\right\},
\label{127b}
\end{equation}
where the quantities $v_{p+}$ appear always in complex conjugate
pairs.
Due to the relations  (\ref{e12a}),(\ref{e12b}),
the polynomial $P_m(u)$ can be written, for even $n$, as a product of factors
of the form  
\begin{eqnarray}
\lefteqn{P_m(u)=\prod_{p=1}^m \left\{e_+(v_+ -v_{p+})+e_-(v_- -v_{p-})
+\sum_{k=1}^{n/2-1}\left\{e_k (v_k-v_{kp})+\tilde e_k(\tilde v_k-
\tilde v_{kp})\right\}\right\}.\nonumber}\\
&&
\label{128a}
\end{eqnarray}
For odd $n$, the polynomial $P_m(u)$ can be written as the product
\begin{eqnarray}
P_m(u)=\prod_{p=1}^m \left\{e_+(v_+ -v_{p+})
+\sum_{k=1}^{(n-1)/2}\left\{e_k (v_k-v_{kp})
+\tilde e_k(\tilde v_k-\tilde v_{kp})\right\}\right\}.
\label{128b}
\end{eqnarray}
These relations can be written with the aid of Eqs. (\ref{e13a}) and
(\ref{e13b}) as
\begin{eqnarray}
P_m(u)=\prod_{p=1}^m (u-u_p) ,
\label{128c}
\end{eqnarray}
where, for even $n$,
\begin{eqnarray}
u_p=e_+ v_{p+}+e_-v_{p-}
+\sum_{k=1}^{n/2-1}\left(e_k v_{kp}+\tilde e_k\tilde v_{kp}\right), 
\label{128d}
\end{eqnarray}
and for odd $n$
\begin{eqnarray}
u_p=e_+ v_{p+}
+\sum_{k=1}^{(n-1)/2}\left(e_k v_{kp}+\tilde e_k\tilde v_{kp}\right), 
\label{128e}
\end{eqnarray}
for $p=1,...,m$.
The roots $v_{p+}$, the roots $v_{p-}$ and, for a given $k$, the roots 
$e_k v_{k1}+\tilde e_k\tilde v_{k1}, ...,  e_k v_{km}+\tilde e_k\tilde v_{km}$
defined in Eqs. (\ref{127a}) or (\ref{127b}) may be ordered arbitrarily.
This means that Eqs. (\ref{128d}) or (\ref{128e}) give sets of $m$ roots
$u_1,...,u_m$ of the polynomial $P_m(u)$, 
corresponding to the various ways in which the roots $v_{p+}, v_{p-}$,
$e_k v_{kp}+\tilde e_k\tilde v_{kp}$ are ordered according to $p$ in each
group. Thus, while the n-complex components in Eq. (\ref{126b}) taken
separately have unique factorizations, the polynomial $P_m(u)$ can be written
in many different ways as a product of linear factors. 

If $P(u)=u^2-1$, the degree is $m=2$, the coefficients of the polynomial are
$a_1=0, a_2=-1$, the n-complex components of $a_2$ are $a_{20}=-1, a_{21}=0,
... ,a_{2,n-1}=0$, the components $A_{2+}, A_{2-}, A_{2k}, \tilde A_{2k}$
calculated according 
to Eqs. (\ref{88a})-(\ref{88d}) are $A_{2+}=-1, A_{2-}=-1, A_{2k}=-1, \tilde
A_{2k}=0, k=1,...,[(n-1)/2]$. The expression of $P(u)$ for even $n$, Eq.
(\ref{126a}), is  
$e_+(v_+^2-1)+e_-(v_-^2-1)+\sum_{k=1}^{n/2-1}\{(e_k v_k+\tilde e_k\tilde
v_k)^2-e_k$\}, and Eq. (\ref{127a}) has the form 
$u^2-1=e_+(v_++1)(v_+-1)+e_-(v_-+1)(v_--1)+
\sum_{k=1}^{n/2-1}\left\{e_k (v_k+1)+\tilde e_k\tilde v_k\right\}
\left\{e_k (v_k-1)+\tilde e_k\tilde v_k\right\}$.
For odd $n$, the expression of $P(u)$, Eq. (\ref{126b}), is  
$e_+(v_+^2-1)+\sum_{k=1}^{(n-1)/2}\{(e_k v_k+\tilde e_k\tilde v_k)^2-e_k\}$,
and Eq. (\ref{127b}) has the form 
$u^2-1=e_+(v_++1)(v_+-1)+
\sum_{k=1}^{(n-1)/2}\left\{e_k (v_k+1)+\tilde e_k\tilde v_k\right\}
\left\{e_k (v_k-1)+\tilde e_k\tilde v_k\right\}$.
The factorization in Eq. (\ref{128c}) is $u^2-1=(u-u_1)(u-u_2)$, where for even
$n$, $u_1=\pm e_+\pm e_-\pm e_1\pm e_2\pm\cdots
\pm e_{n/2-1}, u_2=-u_1$, so that
there are $2^{n/2}$ independent sets of roots $u_1,u_2$
of $u^2-1$. 
It can be checked that 
$(\pm e_+\pm e_-\pm e_1\pm e_2\pm\cdots\pm e_{n/2-1})^2= 
e_++e_-+e_1+e_2+\cdots+e_{n/2-1}=1$.
For odd $n$, $u_1=\pm e_+\pm e_1\pm e_2\pm\cdots
\pm e_{(n-1)/2}, u_2=-u_1$, so that there are $2^{(n-1)/2}$ independent
sets of roots $u_1,u_2$ of $u^2-1$.
It can be checked that 
$(\pm e_+\pm e_1\pm e_2\pm\cdots\pm e_{(n-1)/2})^2= 
e_++e_1+e_2+\cdots+e_{(n-1)/2}=1$.

\section{Representation of polar n-complex numbers by irreducible matrices}

If the unitary matrix written in Eq. (\ref{11}), for even $n$, is called $T_e$,
and the unitary matrix written in Eq. (\ref{12}), for odd $n$, is called $T_o$,
it can be shown that, for even $n$, the matrix $T_e U T_e^{-1}$ has the form 
\begin{equation}
T_e U T_e^{-1}=\left(
\begin{array}{ccccc}
v_+     &     0     &     0   & \cdots  &   0   \\
0       &     v_-   &     0   & \cdots  &   0   \\
0       &     0     &     V_1 & \cdots  &   0   \\
\vdots  &  \vdots   &  \vdots & \cdots  &\vdots \\
0       &     0     &     0   & \cdots  &   V_{n/2-1}\\
\end{array}
\right)
\label{129a}
\end{equation}
and, for odd $n$, the matrix $T_o U T_o^{-1}$ has the form 
\begin{equation}
T_o U T_o^{-1}=\left(
\begin{array}{ccccc}
v_+     &     0     &     0   & \cdots  &   0   \\
0       &     V_1   &     0   & \cdots  &   0   \\
0       &     0     &     V_2 & \cdots  &   0   \\
\vdots  &  \vdots   &  \vdots & \cdots  &\vdots \\
0       &     0     &     0   & \cdots  &   V_{(n-1)/2}\\
\end{array}
\right),
\label{129b}
\end{equation}
where $U$ is the matrix in Eq. (\ref{24b}) used to represent the n-complex
number $u$. In Eqs. (\ref{129a}) and (\ref{129b}), $V_k$ are
the matrices
\begin{equation}
V_k=\left(
\begin{array}{cc}
v_k           &     \tilde v_k   \\
-\tilde v_k   &     v_k          \\
\end{array}\right),
\label{130}
\end{equation}
for $ k=1,...,[(n-1)/2]$, where $v_k, \tilde v_k$ are the variables introduced
in Eqs. (\ref{9a}) and 
(\ref{9b}), and the symbols 0 denote, according to the case, the real number
zero, or one of the matrices
\begin{equation}
\left(
\begin{array}{c}
0  \\
0  \\
\end{array}\right)
\;\; {\rm or}\;\;
\left(
\begin{array}{cc}
0   &  0   \\
0   &  0   \\
\end{array}\right).
\label{131}
\end{equation}
The relations between the variables $v_k, \tilde v_k$ for the multiplication of
n-complex numbers have been written in Eq. (\ref{22}). The matrices 
$T_e U T_e^{-1}$ and $T_o U T_o^{-1}$ provide an irreducible representation
\cite{4} of the n-complex numbers $u$ in terms of matrices with real
coefficients.

\section{Conclusions}

The operations of addition and multiplication of the n-complex numbers
introduced in this 
work have a geometric interpretation based on the amplitude $\rho$,
the modulus $d$ and the polar, planar and azimuthal angles $\theta_+, \theta_-,
\psi_k, \phi_k$. 
If $x_0+x_1+\cdots+x_{n-1}>0$
and $x_0-x_1+\cdots+x_{n-2}-x_{n-1}>0$,
the n-complex numbers can be written in exponential and
trigonometric forms with the aid of the modulus, amplitude and the angular
variables. 
The n-complex functions defined by series of powers are analytic, and 
the partial derivatives of the components of the n-complex functions are
closely related. The integrals of n-complex functions are independent of path
in regions where the functions are regular. The fact that the exponential form
of the n-complex numbers depends on the cyclic variables $\phi_k$
leads to the 
concept of pole and residue for integrals on closed paths. The polynomials of
n-complex variables can be written as products of linear or quadratic
factors.

\newpage

FIGURE CAPTIONS\\

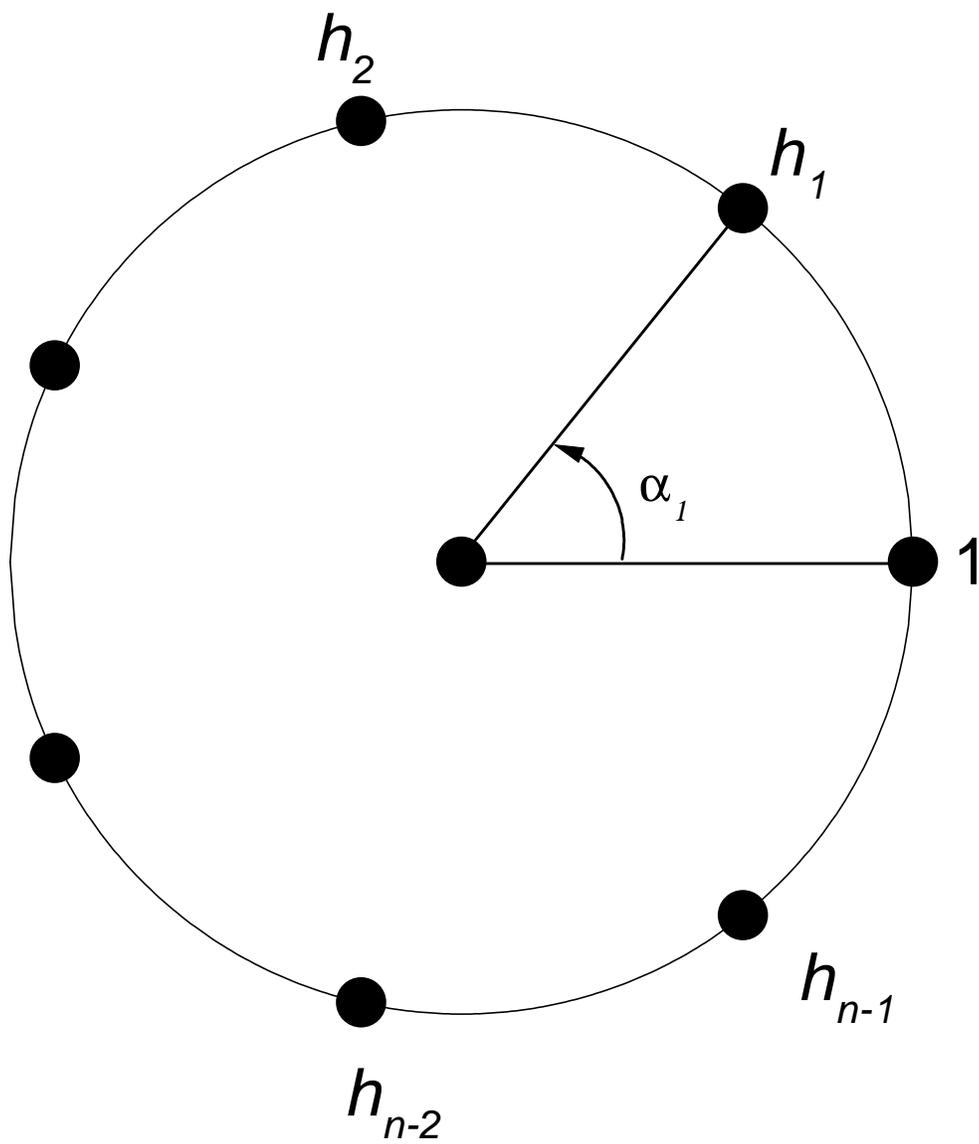
Fig. 1. Representation of the hypercomplex bases $1, h_1,...,h_{n-1}$
by points on a circle at the angles $\alpha_k=2\pi k/n$.
The product $h_j h_k$ will be represented by the point of the circle at the
angle $2\pi (j+k)/n$, $i,k=0,1,...,n-1$. If $2\pi\leq 
2\pi (j+k)/n\leq 4\pi$, the point represents the basis
$h_l$ of angle $\alpha_l=2\pi(j+k)/n-2\pi$.\\

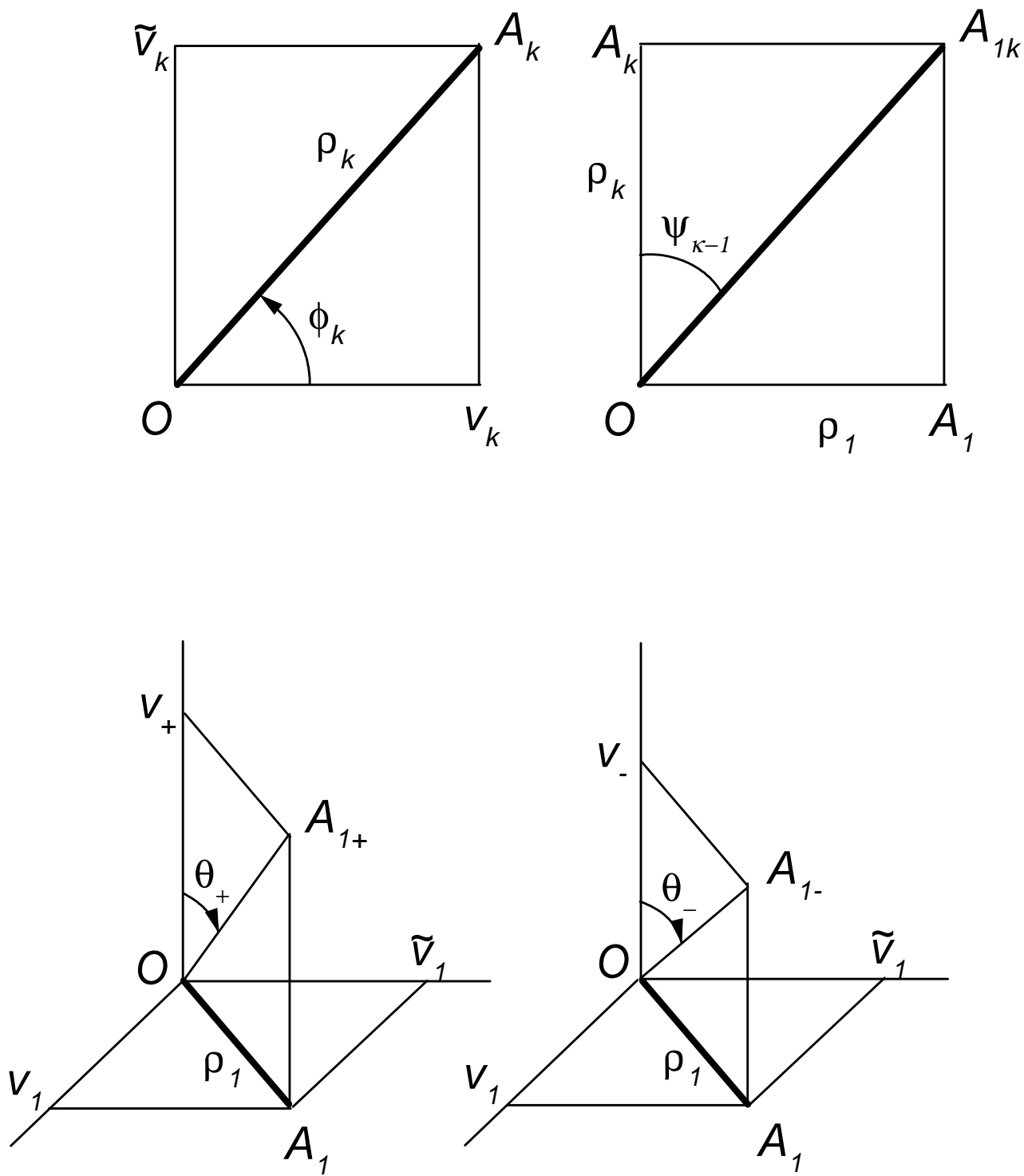
Fig. 2. Radial distance $\rho_k$  and
azimuthal angle $\phi_k$ in the plane of the axes $v_k,\tilde v_k$, and planar
angle $\psi_{k-1}$ between the line $OA_{1k}$ and the 2-dimensional plane
defined by the axes $v_k,\tilde v_k$. $A_k$ is the projection of the point $A$
on the plane of the axes $v_k,\tilde v_k$, and $A_{1k}$ is the projection of
the point $A$ on the 4-dimensional space defined 
by the axes $v_1, \tilde v_1, v_k,\tilde v_k$. 
The polar angle $\theta_+$ is the angle between the line $OA_{1+}$ and the axis
$v_+$, where $A_{1+}$ is the projection of the point $A$ on the 3-dimensional
space generated by the axes $v_1, \tilde v_1, v_+$.
In an even number of dimensions $n$ there is also a polar angle $\theta_-$,
which is the angle between the line $OA_{1-}$ and the axis $v_-$, 
where $A_{1-}$ is the projection of the point $A$ on the 3-dimensional space
generated by the axes $v_1, \tilde v_1, v_-$.\\

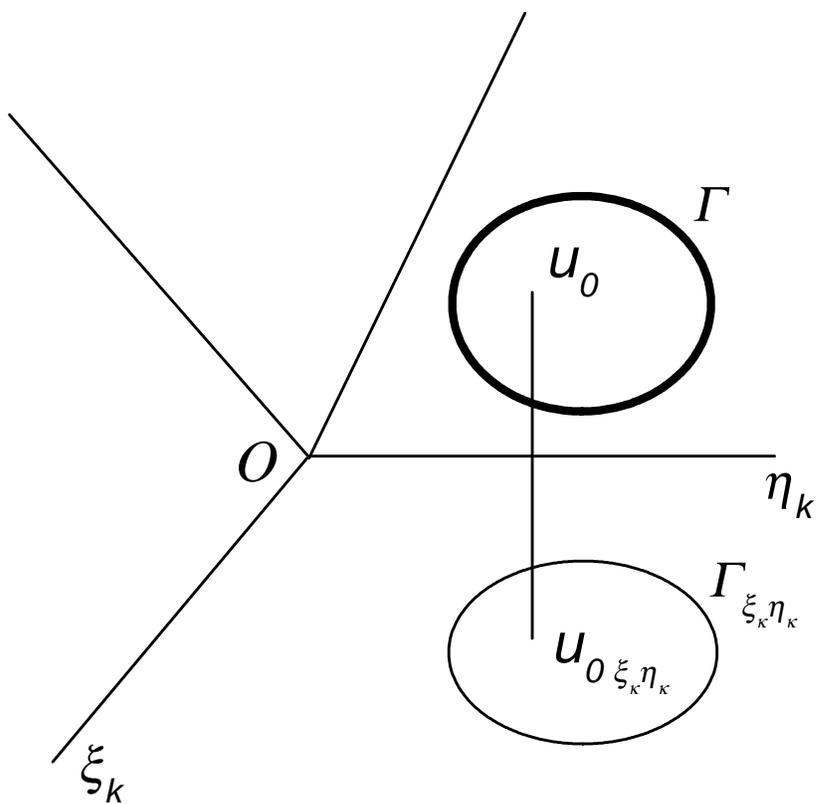
Fig. 3. Integration path $\Gamma$ and pole $u_0$, and their projections
$\Gamma_{\xi_k\eta_k}$ and $u_{0\xi_k\eta_k}$ on the plane $\xi_k \eta_k$.\\ 

\end{document}